\documentclass[12pt]{article}
\usepackage{amsfonts}
\usepackage{amsmath}
\usepackage{amssymb}
\usepackage{latexsym}
\usepackage{amscd}
\usepackage{amsthm}
\usepackage{bbm}
\usepackage{url}

\title{
Adjacency matrices of random digraphs: singularity and anti-concentration}
\author{
Alexander E. Litvak
\and
Anna Lytova
\and
Konstantin Tikhomirov
\and
Nicole Tomczak-Jaegermann
\and
Pierre Youssef
}
\newcommand\address{\noindent\leavevmode

\medskip
\noindent
Alexander E. Litvak, Anna Lytova, Konstantin Tikhomirov
and Nicole Tomczak-Jaegermann,\\
Dept.~of Math.~and Stat.~Sciences,\\
University of Alberta, \\
Edmonton, Alberta, Canada, T6G 2G1.\\
\texttt{\small
e-mail:  aelitvak@gmail.com, lytova@ualberta.ca, ktikhomi@ualberta.ca, \\ nicole.tomczak@ualberta.ca}\\
\medskip

\noindent
Pierre Youssef,\\
Universit\'e Paris Diderot,\\
Laboratoire de Probabilit\'es et de mod\`eles al\'eatoires,\\
75013 Paris, France.\\
\texttt{\small
e-mail:  youssef@math.univ-paris-diderot.fr}
}

\date{}

\newtheorem{theorem}{Theorem}[section]
\newtheorem{lemma}[theorem]{Lemma}
\newtheorem{cor}[theorem]{Corollary}
\newtheorem{defi}[theorem]{Definition}
\newtheorem{prop}[theorem]{Proposition}

\newtheorem{claim}[theorem]{Claim}

\theoremstyle{definition}
\newtheorem{Remark}[theorem]{Remark}

\theoremstyle{plain}

\newtheorem*{theoremA}{Theorem A}
\newtheorem*{theoremB}{Theorem B}

\renewcommand{\proof}{\noindent{\bf Proof.{\ \ }}}

\newcommand{\p}{\mathbb{P}}
\newcommand{\N}{\mathbb{N}}

\newcommand{\R}{\mathbb{R}}

\newcommand{\supp}{{\rm supp\,}}
\newcommand{\rk}{{\rm rk\,}}
\newcommand{\spn}{{\rm span}}

\newcommand{\lam}{\lambda}

\renewcommand{\qed}{\bigskip\hfill\(\Box\)}

\newcommand{\eac}{{\mathcal{E}^{AC}(p)}}

\newcommand{\eoo}{{{\Omega}^{0}}}

\newcommand{\ea}{{\mathcal{E}^{1,2}(q)}}

\newcommand{\erk}{{\mathcal{E}^{1,2}_{rk}}}
\newcommand{\erkij}{{\mathcal{E}^{i,j}_{rk}}}

\newcommand{\ek}{\mathcal{E}_k}
\newcommand{\en}{\mathcal{E}_{n-1}}
\newcommand{\enn}{\mathcal{E}_{n-2}}
\newcommand{\e}{\mathcal{E}}
\newcommand{\K}{{\mathcal{K}^{1,2}(q)}}
\newcommand{\eg}{{\Theta}}

\newcommand{\M}{{M}}
\newcommand{\Mrand}{{M}}

\newcommand{\MMij}{{M}^{i,j}}
\newcommand{\Mc}{\mathcal{M}_{n,d}}
\newcommand{\McUndir}{\mathcal {S}_{n,d}}

\newcommand{\V}{{V}_{1,2}}
\newcommand{\Vij}{{V}_{i,j}}
\newcommand{\Fij}{{F}_{i,j}}

\def\eps{\varepsilon}

\def\UndirIntro{{{\mathcal G}_{n,d}}}
\def\DirIntro{{{\mathcal D}_{n,d}}}
\def\D{{\mathcal{D}}}
\def\Dco{{\mathcal{D}^{co}}}

\def\DO{{\mathcal{D}^0}}
\def\inco{{\rm C_{\nrangr}^{in}}}
\def\outco{{\rm C_{\nrangr}^{out}}}

\def\nrangr{G}
\def\rangr{{G}}
\def\outnbr{{\mathcal N}_{\nrangr}^{out}}
\def\outnbrr{{\mathcal N}_{\nrangr'}^{out}}
\def\outnbrs{{\mathcal N}^{out}}

\def\outnbrpr{{\mathcal N}_{\nrangr'}^{out}}
\def\innbr{{\mathcal N}_{\nrangr}^{in}}
\def\innbrs{{\mathcal N}^{in}}

\def\innbrpr{{\mathcal N}_{\nrangr'}^{in}}
\def\edg{{E}_{\nrangr}}
\def\edgrd{{E}_{\rangr}}
\def\edgpr{{E}_{\nrangr'}}
\def\inedg{{E}_{\nrangr}^{in}}

\def\outedg{{E}_{\nrangr}^{out}}

\def\lvert{\bigl\vert}
\def\rvert{\bigr\vert}

\def\P{{\mathbb P}}

\begin{document}

\maketitle


\begin{abstract}
Let $\DirIntro$ be the set of all $d$-regular directed graphs on $n$ vertices.
 Let $\rangr$ be a graph chosen uniformly at random from $\DirIntro$ and $\Mrand$
be its adjacency matrix. We show that $\Mrand$ is invertible with probability at
least $1-C\ln^{3} d/\sqrt{d}$ for $C\leq d\leq cn/\ln^2 n$, where $c, C$ are
positive absolute constants. To this end, we establish a few properties
of $d$-regular directed graphs. One of them, a Littlewood--Offord type anti-concentration property, is of independent interest.
Let $J$ be a subset of vertices of $\rangr$ with $|J|\approx n/d$.
Let $\delta _i$ be the indicator of the event that the vertex $i$ is connected to
$J$ and define $\delta = (\delta_1, \delta_2, ..., \delta_n)\in \{0, 1\}^n$.  Then for every $v\in\{0,1\}^n$ the probability that $\delta=v$
is exponentially small. This property holds even if a part of the graph is ``frozen."
\end{abstract}

\noindent
{\small \bf AMS 2010 Classification:}
{\small
60C05, 60B20, 05C80, 15B52,
46B06.}

\noindent
{\small \bf Keywords: }
{\small
Adjacency matrices, anti-concentration, invertibility, Littlewood--Offord theory, random digraphs, random graphs, random matrices, regular graphs,
singular probability, singularity, sparse matrices}

\tableofcontents

\section{Introduction}
\label{intro}

For $1\leq d \leq n$ an undirected (resp., directed) graph $G$ is called
{\it $d$-regular} if every vertex has exactly $d$ neighbors (resp., $d$
in-neighbors and $d$ out-neighbors). In this definition we allow graphs to
have loops and, for directed graphs, opposite (anti-parallel) edges, but
no multiple edges. Thus directed graphs ({\it digraphs}) can be viewed as
bipartite graphs with both parts of size $n$. For a digraph $G$ with $n$ vertices
its {\it adjacency matrix} $(\mu_{ij})_{i,j\leq n}$ is defined by
$$
 \mu_{ij}=
 \begin{cases}1,&\mbox{if there is an edge from $i$ to $j$;}\\0,&\mbox{otherwise.}\end{cases}
$$
For an undirected graph $G$  its adjacency matrix  is defined
in a similar way (in the latter case the matrix is symmetric).
We denote the sets of all undirected (resp., directed) $d$-regular graphs
by $\UndirIntro$ and $\DirIntro$, respectively, and the
corresponding sets of adjacency matrices by  $\McUndir$ and $\Mc$.
Clearly $\McUndir \subset \Mc$ and $\Mc$ coincides with the set of $n\times n$ matrices with $0/1$-entries and such that every row and every column  has  exactly $d$ ones. By the probability on $\UndirIntro$, $\DirIntro$, $\McUndir$, and $\Mc$ we always mean  the normalized counting measure.

Spectral properties of adjacency matrices of random $d$-regular graphs
 attracted considerable attention  of researchers
in the recent years.
Among others, we refer the reader to \cite{BHKY}, \cite{BKY}, \cite{DP}, \cite{MR2437174}, \cite{McKay ESD},
 and \cite{TVW} for results
dealing with the eigenvalue distribution.  At the same time, much less
is known about the singular values of the matrices.

The present work is motivated by related general questions on
singular probability.
One problem was mentioned by Vu in his survey \cite[Problem 8.4]{Vu-surv}
(see also 2014 ICM talks by Frieze and  Vu  \cite[Problem~7]{Frieze
  ICM}, \cite[Conjecture~5.8]{Vu}).  It asks if for  $3\leq d\leq n-3$
  the probability that
  a random matrix uniformly distributed on $\McUndir$ is singular goes to zero as
  $n$ grows to infinity.
  Note that in the case $d=1$ the matrix is a permutation matrix, hence non-singular;
  while in the case $d=2$ the conjecture fails (see \cite{Vu-surv} and,  for the directed case, \cite{Cook}). Note also that
  $M\in \Mc$ is singular if and only if the ``complementary" matrix $M'
  \in {\cal{M}}_{n, n-d}$ obtained by interchanging zeros and ones is singular, thus the cases $d=d_0$ and $d=n-d_0$ are essentially the same.
%
The corresponding question for non-symmetric
adjacency matrices is the following (cf., \cite[Conjecture~1.5]{Cook}):

{\it Is it true that for every $3\leq d\leq n-3$
\begin{equation}\tag{$*$}\label{the question}
 p_{n,d}:= \P _{\Mc} \left(\{M\in \Mc:\,M\mbox{ is singular}\} \right)
\longrightarrow
   0\quad \mbox{ as } \quad n\to\infty?
\end{equation}
}

The main difficulty in such singularity questions
stems from the restrictions on row- and column-sums, and from possible
symmetry constraints for the entries.
The question~\eqref{the question} has been recently studied in
\cite{Cook} by Cook who obtained the bound $p_{n,d} \leq d^{-c}$ for
a small universal constant $c>0$ and $d$ satisfying
$\omega(\ln^2 n)\leq d\leq n-\omega(\ln^2 n)$, where $f\geq \omega(a_n)$
means $f/a_n \to \infty$ as $n\to \infty$.

\smallskip

The main result of our paper is the following theorem.

\begin{theoremA}
There are absolute positive constants $c, C$ such that for
 $C\le d\le cn/\ln^2 n$ one has
$$
 p_{n, d} \leq
 \frac{C\ln ^{3}d}{\sqrt{d}}.\;\;
$$
\end{theoremA}

Thus we proved that $p_{n,d} \to 0$ as $d\to \infty$, which in
particular verifies \eqref{the question} whenever $d$ grows to
infinity with $n$, without any restrictions on the rate of convergence.
(Recall that the proof in \cite{Cook} requires $d\geq \omega(\ln^2 n)$.)
We would also like to notice that even for the range
$\omega(\ln^2 n)\leq d\le c n/\ln^2 n$, our bound on probability in Theorem~A is better than in \cite{Cook}.
Of course, it would be nice to obtain a bound going to zero with $n$ and not
with $d$ for the range $d\ge 3$ as well.

\smallskip

In the remaining part of the introduction we describe methods
and techniques used in this paper.  We also explain several novel
ideas that allow us to drop the restriction
$d\geq \omega(\ln^2 n)$ and to treat very sparse matrices. In particular,
we introduce the notion of {\it almost constant} vectors and show how
to eliminate matrices having almost constant null vectors;
we show a new approximation argument dealing with tails of
properly rescaled vectors; we prove an anti-concentration property
for graphs, which is of independent interest; and we provide
a more delicate version of the so-called  ``shuffling" technique.

\smallskip

This paper can be naturally split into two distinct parts.  In the first
one we establish certain properties of random $d$-regular digraphs. In the
second part we use them (or to be more
precise, their ``matrix'' equivalents) to deal with the singularity of
adjacency matrices.  However in the introduction we reverse this
order and discuss first the ``matrix" part as it provides a general
perspective and motivations for graph results.

\smallskip

Singularity of random square matrices is a subject with a long history and many
results. In \cite{Komlos1967} (see also \cite{Komlos Unpublished}) Koml\'os proved
that a random $n\times n$ matrix with independent $\pm 1$ entries
({\it Bernoulli} matrix) is singular
with probability tending to zero as $n\to\infty$.
Upper bounds for the singular  probability of random Bernoulli
matrices were successively improved to $c^n$ (for some $c\in(0,1)$) in
\cite{KKS}; to $\bigl(3/4+o(1)\bigr)^n$ in \cite{TV Bernoulli};  and to
$\bigl(1/\sqrt{2}+o(1)\bigr)^n$ in \cite{BVW}.  Recall that the
conjectured  bound is $\bigl(1/2+o(1)\bigr)^n$. The corresponding problem for symmetric Bernoulli matrices was considered
in \cite{CTV}, \cite{Nguyen}, \cite{Vershynin Symmetric}.
Recently, matrices with independent rows and with row-sums constrained
to be equal to zero were studied in \cite{Nguyen Constant sum}.

In all these works, a fundamental role is played by what is nowadays
called {\it the Littlewood-Offord theory}. In its classical form,
established by Erd{\H{o}}s \cite{Er}, the Littlewood-Offord inequality states that for every fixed $z\in\R$, a fixed vector $a=(a_1,a_2,\dots,a_n)\in\R^n$
with non-zero coordinates, and for independent random signs $r_k$, $k\le n$, the probability $\P\left\{\sum_{k=1}^n r_ka_k=z\right\}$
is bounded from above by $n^{-1/2}$.
This combinatorial result has been substantially strengthened and
generalized in subsequent years, leading to a much better
understanding of interrelationship between
the law of the sum $\sum_{k=1}^n r_ka_k$ and the arithmetic
structure of the vector $a$.
For more information and further references, we refer the
reader to \cite{TaoVu}, \cite[Section~3]{TV LO to Circular}, and \cite[Section~4]{RV  Congress}.
The use of the Littlewood-Offord theory in context of random matrices
can be illustrated as follows. Given an $n\times n$ matrix $A$ with
i.i.d.\ elements, $A$ is non-singular if and only if the inner product
of a normal vector to the span of any subset of $n-1$ columns of $A$
with the remaining column is non-zero. Thus, knowing the ``typical''
arithmetic structure of the random normal vectors and conditioning on
their realization, one can estimate the probability that $A$ is
singular. Moreover, a variant of this approach allows us to obtain sharp
quantitative estimates for the smallest singular value of the matrix
with independent subgaussian entries \cite{RV Square}.

Similarly to the aforementioned works, the Littlewood-Offord
inequality plays a crucial role in the proof of Theorem~A.
Note that if $\Mrand$ is a random matrix uniformly distributed on
$\Mc$ then every
two entries/rows/columns of $\Mrand$ are probabilistically dependent;
moreover, a realization of the first $n-1$ columns uniquely defines
the last column of $\Mrand$.  This makes a straightforward application
of the Littlewood-Offord
theory (as illustrated in the previous paragraph) impossible.

In \cite{Cook}, a sophisticated approach based on the ``shuffling'' of two
rows was developed to deal with that problem.
The shuffling consists in  a random perturbation of two rows of a
fixed matrix $M\in\Mc$ in such a way that the sum of the rows remains
unchanged. We discuss this procedure in more details in
Section~\ref{main-result}. It can be also defined in terms of ``switching"
discussed below.
The proof in \cite{Cook} can be divided into two steps: at the first
step, one proves that the event  that a random matrix $\Mrand$
does not have any (left or right)
null vectors with many ($\geq C n d^{-c}$) equal coordinates has
probability close to one, provided that $d\ge \omega (\ln^2 n)$.
Then one shows that conditioned on this event, a random matrix $\Mrand$ is
non-singular with large probability.

In our paper, we expand on  some of the techniques developed in
\cite{Cook} by adding new crucial ingredients.
On the first step, in Section~\ref{s:AC}, we show
that for $C\leq d\leq cn$, with probability going to one with $n$, a random
matrix $\Mrand$ does not have any null vectors having at least
$n(1-1/\ln d)$  equal coordinates,  (we call such vectors {\it almost constant}). Note that we rule out a much smaller set of null vectors. This  allows us to drop the lower bound on $d$, but requires a delicate
adjustment of the second step.
Key elements of the first step consists of a new {\it anti-concentration}
property of random graphs and their adjacency matrices as well as of
using a special form of an $\varepsilon$-net
build  from  the ``tails'' of appropriately rescaled
vectors $x\in\R^n$.
Then, conditioning on the event that $\Mrand$ does not have
almost constant null vectors, we show in Section~\ref{main-result}
that  a random matrix $\Mrand$ is non-singular with high probability.
This relies on a somewhat modified and simplified version of the
shuffling procedure for the matrix rows.
As the shuffling involves supports of only two rows we get at this step that
probability converges with $d$ and not with $n$. We would like
to emphasize that this is the only step which does not allow to have the  convergence to zero with $n$.


\smallskip

We now turn our attention to Section~\ref{graphs}, which deals with the set  $\DirIntro$ of $d$-regular digraphs.
Our analysis
is based on an operation called ``the simple switching,''
which is a standard tool to work with regular graphs.  As an
illustration, let $G\in\DirIntro$ and let $i_1\neq i_2$ and
$j_1\neq j_2$ be vertices of $G$ such that $(i_1,j_1)$ and $(i_2,j_2)$
are edges of $G$ and $(i_1,j_2)$, $(i_2,j_1)$ are not.  Then the
simple switching consists in replacing the edges $(i_1,j_1)$,
$(i_2,j_2)$ with $(i_1,j_2)$ and $(i_2,j_1)$, while leaving all other
edges unchanged. Note that the operation does not destroy
$d$-regularity of the graph.  The simple switching was introduced (for
general graphs) by Senior \cite{Senior} (in that paper, it was called
``transfusion''); in the context of $d$-regular graphs it was first
applied by McKay \cite{McKay subgraphs}. As in \cite{McKay subgraphs},
we use this operation to compare cardinalities of certain subsets of $\DirIntro$.
We note that one could use the configuration model, introduced by Bollob\`as \cite{MR595929} in the context of random regular graphs,
to prove our results for sparse graphs. We prefer to use the switching method in order to have a unified proof
for all ranges of $d$.

As in the matrix counterpart we work with a random graph $\rangr$
uniformly distributed on $\DirIntro$.
For a finite set $S$, we denote by $|S|$ its cardinality.
For a positive integer $n$ we denote by $[n]$ the set $\{1,2,\dots,n\}$.
For every subset $S\subset [n]$, let $\innbrs_\rangr(S)$ be the set of all vertices
of $\rangr$ which are in-neighbors to some vertex in $S$. Further, for every two
subsets $I,J$ of $[n]$, we denote by $\edgrd(I,J)$ the collection of edges of $\rangr$ starting from a vertex in $I$ and ending at a vertex from $J$.
In a simplified form, our first statement about graphs
(Theorem~\ref{th-graph-sparse} in
Section~\ref{section-isoper-index-graph})
can be formulated as follows:

\smallskip

{\it Let $8\leq d\leq n$, $\varepsilon\in (0,1)$, and $k\geq 2$. Assume that
$\varepsilon^2 \geq d^{-1} \max\{8, \ln d\}$ and $k\leq c\varepsilon n/d$
for a sufficiently small absolute positive constant $c$.
Then}
$$
   \P\bigl\{\exists S\subseteq [n],  \,
\vert S\vert =k \, \text{ such that }\, \,
  \lvert\innbrs_\rangr(S)\rvert  \leq  (1-\varepsilon)d\lvert S\rvert \bigr\} \leq
    \exp\left(-\frac{\eps^2 d k}{8}\, \ln\left(\frac{3ec\eps n}{k d}\right)\right).
$$

Note that $|S|\leq \lvert\innbrs_\rangr(S)\rvert\leq d|S|$.
Thus, roughly speaking, our result says
 that ``typically,'' whenever a set $S$ is not too large,
the set of all in-neighbors of $S$ has cardinality close to the maximal possible one.
In the
case of undirected graphs such results are known (see e.g. \cite{Noga-book} and references therein).
We note that in fact we prove a more general statement, in which we estimate the probability
conditioning on a ``partial'' realization of a random graph $\rangr$, when a certain subset of its edges is fixed
(see Theorem~\ref{th-graph-sparse}).

In our second result, we estimate the probability that $\edgrd(I,J)$ is empty for
large sets $I$ and $J$ (see Theorem~\ref{th-graph-spread} in
Section~\ref{section-independence-graph}):

\smallskip

{\it
There exist absolute positive constants $c,  C$
such that the following holds. Let $2\leq d\leq  n/24$ and  $C n \ln d/d  \leq \ell \leq  r\leq n/4$. Then}
$$
 \P\bigl\{\edgrd(I,J)=\emptyset\;\;\mbox{for some }I,J\subset[n]\mbox{ with }|I|\geq \ell,\;|J|\geq r\bigr\}
 \leq \exp\left(-c r\ell d/n\right).
$$

Note that the first statement can be  reformulated in terms of sets
$\edgrd(I,J)$ (however, the range of cardinalities for $I$ and $J$
will be different compared to the second result). These statements
can be seen as manifestations of a general phenomenon that a random
graph $\rangr$ with a large probability has good regularity properties.
Let us also note that analogous statements for the Erd{\H{o}}s--R\'enyi graphs
(in this random model an edge between every two vertices is included/excluded
in a graph independently of other edges) follow from standard Bernstein-type inequalities. For related results on $d$-regular random graphs, we refer
the reader to \cite{MR1839497} where concentration properties of {\it co-degrees}
were established in the undirected setting, and to \cite{Cook Graphs} for concentration of co-degrees and of
the ``edge counts'' $|\edgrd(I,J)|$ for digraphs. In paper \cite{Cook Graphs} which
serves as a basis for the main theorem of \cite{Cook} mentioned above, rather strong concentration properties of $|\edgrd(I,J)|$ are established;
however, the results provided in that paper are valid only for
$d\geq \omega(\ln n)$.
The proof in \cite{Cook} is based on the method of exchangeable pairs introduced by Stein and developed
for concentration inequalities by Chatterjee (see survey \cite{Chatterjee} for more information and references).
On the contrary, our proof of the afore-mentioned statements is simpler,  completely self-contained and works for $d\geq C$.
As we mentioned above, we use the following
Littlewood-Offord type anti-concentration result
matching anti-concentration properties of a
weighted sum of independent random variables or vectors
studied in the Littlewood-Offord theory.
This result is of independent interest, and we
formulate it here as a theorem (see also Theorem~\ref{th-anticoncentration} in
Section~\ref{section-anticoncentration-graph}).
For every $J\subset [n]$ and $i\in[n]$ we define
$\delta_i^{J}(\rangr)\in \{0,1\}$ as the indicator
of the event $\{i\in \innbrs_\rangr(J)\}$ and denote
$\delta^J(\rangr):=(\delta_1^J(\rangr),\ldots, \delta_n^J(\rangr))\in
\{0,1\}^n$.
\begin{theoremB}
There are two positive absolute constants $c$ and $c_1$ such that the
following holds.
Let $32\leq d\leq cn$ and  $I, J$ be disjoint subsets of $[n]$ such that
$\vert I\vert \leq d\vert J\vert/32$ and $8\leq \vert J \vert  \leq 8 c n/d$.
 Let vectors $a^i\in\{0,1\}^n$, $i\in I$, be
such that the event
$$
 \e:=\{\innbrs_\rangr(i)=\supp a^i\mbox{ for all }i\in I\}
$$
has non-zero probability (if $I=\emptyset$ we set $\e=\DirIntro$).
Then for every  $v\in \{0,1\}^n$ one has
$$
  \P\{ \delta^J(\rangr) =v  \mid  \e\}\leq 2\exp\left(- c_1 d
  \vert J\vert \ln\big(\frac{n}{ d\vert J\vert }\big)\right).
$$
\end{theoremB}

We note that the probability estimate in the previous statement
matches the one for the corresponding quantity $\delta^J$ in the Erd{\H{o}}s--R\'enyi model.

\medskip

The paper is organized as follows.
Sections~\ref{graphs} deals with all results related to graphs.
 Section~\ref{fromgraph} provides links  between the graph
results of Section~\ref{graphs} and the matrix results used in
Section~\ref{s:matrices}.
Finally, Section~\ref{s:matrices} presents the proof of the main theorem,
including a  number of auxiliary combinatorial lemmas.

\smallskip

In this paper letters $c, C$, $c_0, C_0$, $c_1, C_1$, ... always denote
absolute positive constants (i.e. independent of any parameters), whose
precise value may be different from line to line.

\smallskip

Main results of this paper were announced in \cite{LLTTY}.

\bigskip

\noindent
{\bf Aknowledgment. } This work was conducted while the second named author was a
    Research     Associate at the  University of Alberta,
    the third named author was a graduate student and held
        the PIMS Graduate Scholarship,
    and the last named author was a CNRS/PIMS PDF at
             the same  university.
    They all  would like  to thank the  Pacific Institute
        and the University of Alberta for the support.
        A part of this work was also done when the first four authors
took part in activities of the annual program ``On Discrete Structures: Analysis and Applications" at the Institute for Mathematics and its Applications (IMA), Minneapolis, MN, USA.
These authors would like to thank IMA
for the support and excellent working conditions.
    All authors would like to thank Michael Krivelevich for  many
    helpful comments  on the ``graph" part of this paper. We  would also like to thank Justin Salez for helpful comments.

\section{Expansion and anti-concentration for random digraphs}
\label{graphs}


\subsection{Notation and preliminaries}


For a real number $x$, we  denote by $\lfloor x\rfloor$ the largest integer smaller than or equal to $x$ and
by $\lceil x\rceil$ the smallest integer larger than or equal to $x$.
Further, for every $a\geq 1$, we denote by $[a]$ the set
$\{i\in\N:\, 1\leq i\leq \lfloor a\rfloor\}$.

Let $d\leq n$ be positive integers. A {\it $d$-regular directed graph} (or
{\it $d$-regular digraph}) on $n$ (labeled) vertices is a graph in which every
vertex has exactly $d$ in-neighbors and $d$ out-neighbors. We allow the graphs
to have loops and opposite/anti-parallel edges but do not allow multiple edges.
Thus this set coincides with the set of $d$-regular bipartite graphs with both
parts of size $n$. The set of vertices of such graphs is always  identified with $[n]$.
The set of all these graphs is denoted by $\DirIntro$. When $n$ and $d$ are clear
from the context, we will use a one-letter notation $\D$. Note that the set of adjacency matrices for graphs in $\D$ coincides with the set of $n\times n$
matrices with $0/1$-entries such that every row and every column has exactly $d$ ones.
By a random graph on $\D$ we always
mean a graph uniformly distributed on $\D$ (that is, with respect to the normalized counting measure).

Let $\nrangr=([n],E)$ be an element of $\D$, where $E$ is the set of its
directed edges.
Thus $(i,j)\in E$, $i,j\leq n$, means that there is an edge going from vertex $i$ to vertex $j$.
We will denote the adjacency
matrix of $\nrangr$ by $M=M(\nrangr)$; its rows and columns by $R_i=R_i(M)=R_i(\nrangr)$ and
$X_i=X_i(M)=X_i(\nrangr)$, $i\leq n$, respectively.

\smallskip


Given a graph $\nrangr\in\D$ and a subset $S\subset [n]$ of its vertices, let
\begin{align*}
\outnbrs(S)&=\outnbr(S):=\bigl\{v\leq n:\,  \exists i\in S \, \, (i,v) \in E \bigr\} = \bigcup _{i\in S} {\rm supp} R_i,\\ 
\innbrs(S)&=\innbr(S):=\bigl\{v\leq n:\,  \exists i\in S \, \, (v,i) \in E
\bigr\}=\bigcup _{i\in S} {\rm supp} X_i.
\end{align*}
Similarly, we define the out-edges and the in-edges as follows
\begin{align*}
\outedg(S)&:=\bigl\{e\in E:\, e=(i,j) \text{ for some } i\in S \text{ and } j\leq n \bigr\},\\
\inedg(S)&:=\bigl\{e\in E:\, e=(i,j) \text{ for some } i\leq n \text{ and } j\in S\bigr\}.
\end{align*}
For one-element subsets of $[n]$ we will use lighter notations
$\outnbr(i)$,
$\innbr(i)$, $\outedg(i)$, $\inedg(i)$ instead of
$\outnbr(\{i\})$,
$\innbr(\{i\})$, $\outedg(\{i\})$, $\inedg(\{i\})$.

\smallskip

Given a graph $\nrangr=([n],E)$, for every $I,J\subset [n]$ the set of all edges departing from $I$ and landing in $J$ is denoted by
$$
\edg(I\times J)= \edg(I,J)=\bigl\{ e\in E:\, e=(i,j) \text{ for some } i\in I \text{ and } j\in J\bigr\}.
$$
Further, we let
$$
\DO(I,J) = \bigl\{ \nrangr\in \D:\, \edg(I, J)=\emptyset\bigr\}.
$$
Note that $\DO(I,J)$ is the set of all graphs whose adjacency matrices have zero $I\times J$-minor,
hence the superscript ``$0$''.

Given $\nrangr\in\D$, for $u, v\leq n$ the sets of common out-neighbors and common in-neighbors will be denoted as
\begin{align*}
\outco(u, v)&=\{j\leq n\, :\, (u, j), (v, j)\in E\} = \supp R_u \cap \supp R_v,\\
\inco(u, v)&=\{i\leq n \, :\, (i, u), (i, v)\in E\} = \supp X_u \cap \supp X_v.
\end{align*}

For every $S\subset [n]$ and $F\subset [n]\times [n]$, we define
$$
\D(S,F)=\bigl\{ \nrangr\in\D:\, \inedg( S)=F\bigr\}.
$$
Informally speaking, $\D(S,F)$ is the subset of $d$-regular graphs for which
the in-edges of $S$ are ``frozen'' and, as a set, coincide with $F$.
Note that a necessary (but not sufficient) condition for $\D(S,F)$ to be non-empty is
$$
\forall i\leq n\, \, \, \, \vert \{\ell\in [n]:\, (i, \ell)\in F\}\vert \leq d \quad \mbox{ and } \quad  \forall j\in S\, \,  \, \,
\vert \{\ell\in [n]:\, (\ell, j)\in F\}\vert = d.
$$

For every $\varepsilon \in (0,1)$, denote
$$
  \Dco(\varepsilon)=\bigl\{ \nrangr\in\D:\, \forall i\neq j\leq n\,\,\, \, \, \,
  \vert \outco(i,j)\vert\leq \varepsilon d\bigr\} = \bigcap_{i<j}
   \D_{i,j}^{co}(\varepsilon),
$$
where
$$
   \D_{i,j}^{co}(\varepsilon):=\bigl\{ \nrangr\in\D:\, \vert \outco(i,j)\vert\leq
   \varepsilon d\bigr\}.
$$

\bigskip

Let $A$, $B$ be sets, and $R\subset A\times B$ be a relation.
Given $a\in A$ and $b\in B$, the image of $a$ and preimage of $b$ are defined by
$$
     R(a) = \{ y \in B \, : \, (a,y)  \in  R\} \quad \mbox{ and } \quad
     R^{-1}(b) = \{ x \in A \, : \, (x,b)  \in  R\}.
$$
We also set $R(A)=\cup _{a\in A} R(a)$.
Further in this section, we often define relations between sets in order to estimate their cardinality,
using the following simple claim.

\begin{claim} \label{multi-al}
Let $s, t >0$.
Let $R$ be a relation between two finite sets $A$ and $B$ such that for
every $a\in A$ and every $b\in B$ one has $|R(a)|\geq s$ and $|R^{-1}(b)|\leq t$.
Then
$$
        s |A|\leq  t |B|.
$$
\end{claim}

\proof
 Without loss of generality we assume that $A=[k]$ and $B=[m]$ for some positive integers $k$ and $m$.
For $i\leq k$ and
$j\leq m$, we set $r_{ij} =1$ if $(i,j)\in R$ and $r_{ij} =0$ otherwise. Counting the number
of ones in every row and every column of the matrix $\{r_{ij}\}_{ij}$ we obtain
$$
      \sum _{i=1}^k \sum _{j=1}^m r_{ij} =  \sum _{i=1}^k |R(i)| \geq  s k = s |A|
      \quad \mbox{ and } \quad
      \sum _{j=1}^m \sum _{i=1}^k r_{ij} = \sum _{j=1}^m |R^{-1}(j)| \leq t m = t |B|,
$$
which implies the desired estimate.
\qed


\subsection{An expansion property of random  digraphs}
\label{section-isoper-index-graph}

In this section, we
establish certain expansion properties of random graphs uniformly distributed on $\D$, which
can roughly be described as follows: given a subset $S\subset [n]$ of cardinality $|S|\leq c n/d$,
with high probability the number of in-neighbors of $S$ is of order $d|S|$.
Beside its own interest, this result is used in the proof of the anti-concentration property for graphs which will be given in
Section~\ref{section-anticoncentration-graph}. In fact we will need a
statement where we control the number of in-neighbors of
a subset of vertices while ``freezing'' (i.e.\ conditioning on a realization of) a set of edges inside the graph.

\begin{theorem}\label{th-graph-sparse} Let $8\leq d\leq n$, $\varepsilon\in (0,1)$, and $k\geq 2$. Assume that
$$
\varepsilon^2 \geq \frac{\max\{8, \ln d\}}{d}\quad\text{and}\quad k\leq \frac{c\varepsilon n}{d}
$$
for a sufficiently small absolute positive constant $c$.
Let $I\subset [n]$ be of cardinality at most $n/8$.
Define
$$
  \Gamma _k=\bigl\{ \nrangr\in\D \, :\, \exists S\subseteq I^c,  \,
\vert S\vert =k, \, \, \text{ such that }\, \,
  \lvert\innbr(S)\rvert  \leq  (1-\varepsilon)d\lvert S\rvert  \bigr\}
$$
and
$$
  \Gamma=\bigl\{ \nrangr\in\D\, :\, \exists S\subseteq I^c, \,
  \vert S\vert \leq c\eps n/d,\, \,  \text{ such that }\, \,
  \lvert\innbr(S)\rvert \leq (1-\varepsilon)d\lvert S\rvert  \bigr\}
  = \bigcup _{\ell=2}^{c\eps n/d} \Gamma _\ell .
$$
Then for every $F\subset [n]\times [n]$ with $\D(I,F)\neq \emptyset$ we have
$$
   \P\left(\Gamma _k\mid \D(I,F) \right) \leq
    \exp\left(-\frac{\eps^2 d k}{8}\, \ln\left(\frac{3ec\eps n}{k d}\right)\right) .
$$
In particular,
$$
   \P\left(\Gamma \mid \D(I,F) \right) \leq
   \exp\left(-\frac{\eps^2 d}{8} \ln\left(\frac{ec\eps n}{d}\right)\right) .
$$
\end{theorem}

Let us describe the idea of the proof of Theorem~\ref{th-graph-sparse}. Suppose we are given
a set of vertices $S$ of an appropriate size. Since $\vert \inedg(S)\vert =d\vert S\vert$, then
we always have
$$
\vert S\vert \leq \vert \innbr(S)\vert\leq d\vert S\vert.
$$
We want to prove that the number of graphs satisfying $\vert \innbr(S)\vert\leq (1-\varepsilon)d\vert S\vert$
is rather small.
In order to estimate the number of in-neighbors of $S$, our strategy is to build
$S$ by adding one vertex at a time and trace how the number of
in-neighbors is changing. Namely, if $S=\{v_i\}_{i\leq s}$ then to build $S$ we start by setting $S_1:=\{v_1\}$ -- a set for which we know that it has exactly $d$ in-neighbors. Now we add the vertex
$v_2$ to $S_1$ to get $S_2:=\{v_1,v_2\}$. We need to trace how
the number of in-neighbors to $S_2$ changed compared to that of $S_1$. More precisely,
we need to count the number of graphs for which the number of in-neighbors
has increased by at most $(1-\varepsilon/2)d$. To this end,
we count the number of graphs having the property that the number of common in-neighbors to $v_1$ and $v_2$ is at least
$\varepsilon d/2$. We count such graphs by applying the simple switching.
One should be careful here to switch the edges without interfering with the frozen area of the graph.
We continue in a similar manner by adding one vertex at a time and controlling the number of
common in-neighbors between the added vertex and the existing ones.
Now, note that the condition $\vert \innbr(S)\vert \leq (1-\varepsilon)d\vert S\vert$
implies that for a large proportion of the vertices added, the number of common in-neighbors with the existing vertices is
at least $\varepsilon d/2$. We use this together with the cardinality estimates obtained via the simple switching at
each step to get the required result.

We  use the following notation.
Given $S\subset [n]$ and $\delta\in (0,1)$,
we set $\Gamma(S, \emptyset)=\D$ and for a non-empty $J\subset [n]$, let
$$
   \Gamma(S, J)=\Gamma(S, J, \delta)=\Bigl\{G\in\D \, : \, \forall j\in J \, \mbox{ one has } \Big| \bigcup_{i\in S, i<j} \inco(i, j) \Big| \geq \delta d\Bigr\}
$$
(the number $\delta$ will always be clear from the context). We  also use a simplified notation $\Gamma(S, j):= \Gamma(S, \{j\})$.

Note that $\Gamma(S, J)$ contains all graphs in which every vertex $j\in J$ has many common in-neighbors
with the set $\{i\in S\, \, :\, \, i<j\}$. In the next lemma, we  estimate cardinalities of $\Gamma(S,j)$, conditioning
on a ``partial'' realization of a graph.

\begin{lemma}\label{lem2-graph-sparse}
Let $\delta \in (0,1)$, $2  \leq d\leq  n/12$,
$1\leq k\leq  \delta n/{(4e d)}$ and
 $F, H\subset [n]\times[n]$. For every $I\subset [k+d]^c$ satisfying
 $$
 \vert I\vert \leq \frac{n}{8},
 $$
one has
$$
  \lvert \Gamma([k], k+1)\cap \D([k],F)\cap \D(I,H) \rvert
 \leq \gamma_k \lvert \D([k],F)\cap \D(I,H)\rvert ,
$$
where
$$
 \gamma _k= \left(\frac{2 ekd}{\delta n}\right)^{\delta d}.
$$
\end{lemma}
Less formally, the above statement asserts that, considering a subset of $\D$
with prescribed (frozen) sets of in-edges for $[k]$ and $I$, for a vast majority of such graphs
the $(k+1)$-th vertex will have a small number of common in-neighbors with the first $k$ vertices.

\medskip

\proof
We assume that the intersection $\D([k],F)\cap \D(I,H)$ is non-empty.
Then we have $F([n])=[k]$ and $F^{-1}([k])=\innbr ([k])$
(recall notation for images and preimages of a relation).
 Without loss of generality, $\innbr ([k]) = [n_1]^c$ for some $n_1\leq n$.
Note that
$$
  k\leq \big\vert \innbr ([k])\big\vert \leq kd,
$$
hence $n-kd \leq n_1\leq  n-k$.

For $0\leq q\leq d$ denote
$$
 Q(q):=\bigl\{ \nrangr \in \D([k],F)\cap \D(I,H):\, \big\vert
 \innbr([k])\cap \innbr(k+1) \big\vert =q\bigr\}.
$$
and
$$
  Q:= \Gamma([k], k+1) \cap \D([k],F)\cap \D(I,H)=\bigcup _{q=\lceil\delta d\rceil}^d Q(q).
$$

\smallskip

We proceed by comparing the cardinalities $Q(q)$ and $Q(q-1)$ for every $1\leq q\leq d$.
To this end, we will define a relation $R_q\subset Q(q)\times Q(q-1)$.
Let $\nrangr \in Q(q)$. Then there exist
$n_1 <i_1<...<i_{q}$ such that for every $\ell\leq q$ we have
$$
 i_{\ell}\in \innbr([k])\cap \innbr(k+1).
$$

For every $\ell \leq q$,  there are at most $d^2$ edges inside
$
\edg\big([n_1],  \outnbr(i_\ell)\big).
$
Further, there are $(n_1-(d-q))d$ edges in $\outedg\big([n_1]\setminus \innbr(k+1)\big)$ and
at most $d\vert I\vert$ edges in $\edg\big([n_1]\setminus \innbr(k+1), I\big)$.
Therefore, for every $\ell \leq q$, the cardinality of the set
$$
E_\ell:=\edg\big([n_1]\setminus \innbr(k+1), I^c\setminus \outnbr(i_\ell)\big)
$$
can be estimated as
$$
  |E_\ell|\geq \big(n_1-(d-q)-\vert I\vert \big)d-d^2\geq (7n/8 - kd - 2d) d\geq  nd/2
$$
(here, we used the conditions $|I|\leq n/8$ and $n_1\geq n-kd$ together with the restrictions on $k$).

Now, we turn to constructing the relation $R_q$.
We let a pair $(\nrangr, \nrangr')$ belong to $R_q$ for some $\nrangr' \in Q(q-1)$ if
$\nrangr'$ can be obtained
from $\nrangr$ in the following way. First we choose $\ell \leq q$ and an edge $(i,j)\in E_\ell$. We destroy the edge
$(i_\ell, k+1)$ to form the edge $(i,k+1)$, then destroy the edge $(i,j)$ to form the edge $(i_\ell,j)$
(in other words, we perform the simple switching on the vertices $i,i_\ell,j,k+1$). Note that
the conditions $i\notin \innbr(k+1)$ and $j\notin \outnbr(i_\ell)$, which are implied by the definition of $E_\ell$,
guarantee that the simple switching does not create multiple edges, and we obtain a valid graph in $Q(q-1)$.

The definition of $R_q$ implies that for every $G\in Q(q)$ one has
\begin{equation}\label{eq-image-Rq-lem51}
\vert R_q(\nrangr)\vert\geq\sum\limits_{\ell=1}^q |E_\ell|\geq \frac{qnd}{2}.
\end{equation}

Now we estimate the cardinalities of  preimages.
Let $\nrangr' \in R_q\bigl( Q(q)\bigr)$.
In order to reconstruct a graph $\nrangr$ for which
$(\nrangr, \nrangr')\in R_q$, we need to perform a simple switching which
$$
  \mbox{destroys an edge from} \quad \edgpr \bigl([n_1], k+1\bigr) \quad \quad
  \mbox{and}\quad \quad \mbox{adds an edge to} \quad \edgpr \bigl([n_1]^c, k+1\bigr)  .
$$
There are at most $d-q+1$ choices to destroy an edge in $\edgpr \bigl([n_1], k+1\bigr)$,
and at most $n-n_1\leq kd$ possibilities to create an edge connecting $[n_1]^c$ with $(k+1)$-st vertex.
Assume that we destroyed an edge $(v, k+1)$ and added an edge $(u, k+1)$.
The second part of the simple switching is to destroy an excessive out-edge
of $u$ and create a corresponding edge (with the same end-point) for $v$. It
is easy to see that we have at most $d$ possibilities for the second part
of the switching. Therefore,
\begin{equation*}
\lvert R_q^{-1}(\nrangr')\rvert \leq kd^3.
\end{equation*}

Using this bound, Claim~\ref{multi-al}, and \eqref{eq-image-Rq-lem51},
we obtain that
$$
  \lvert Q(q)\rvert \leq \left(\frac{2 kd^2}{qn}\right)\cdot  \lvert Q(q-1)\rvert
$$
and, applying the estimate successively,
$$
    \lvert Q(q)\rvert     \leq          \left(\frac{2 kd^2}{n}\right)^{q} \, \,
    \frac{1}{q!}\, \,   \lvert Q(0)\rvert.
$$
Since $q!\geq 2(q/e)^q$ and $2e kd/(\delta n)\leq 1/2$, this implies
$$
\vert Q\vert= \sum_{q=\lceil\delta d\rceil}^d\lvert Q(q)\rvert
\leq \frac{1}{2}\,\,  \sum_{q=\lceil\delta d\rceil}^d \left(\frac{2e kd}{\delta n}\right)^{q} \,
\lvert Q(0)\rvert
\leq
 \left(\frac{2e kd}{\delta n}\right)^{\delta d}\,
\lvert Q(0)\rvert .
$$
Using that $Q(0)\subset \D([k],F)\cap \D(I,H)$, we obtain the desired result.
\qed

Now, we iterate the last lemma to obtain the following statement.

\begin{cor}\label{cor-al}
Let $\delta$, $n$,  $d$, $k$ and  $\gamma _k$ be as in Lemma~\ref{lem2-graph-sparse} and let $\ell\leq k$.
Further, let $I\subset [n]$ satisfy
$
\vert I\vert \leq n/8
$
and let $H\subset [n]\times[n]$.
Then for every subsets $J\subset S\subset I^c$ such that $\vert S\vert =k$ and $\vert J\vert =\ell$, one has
$$
 \lvert \Gamma(S, J)\cap \D(I,H) \rvert \leq  {\gamma_k}^{\ell} \,  \lvert \D(I,H) \rvert .
$$
\end{cor}

\proof
Without loss of generality we assume that the intersection
$\Gamma(S, J)\cap \D(I,H)$ is non-empty, that $S=[k]$ and $I\subset [k+d]^c$.
Write $J=\{j_1, ..., j_{\ell}\}$ for some $j_1<...<j_{\ell}$. For $1\leq s\leq \ell$ denote
$J_s=\{j_1, ..., j_{s}\}$, $J_0=\emptyset$ and let $k_s=j_{s}-1$.
Note that for every $1\leq s\leq \ell$, we have
$$
\Gamma (S, J_s)= \Gamma ([k_s], J_s).
$$
Note also that
\begin{equation}\label{eq-coro-split}
\Gamma (S, J_s)= \Gamma ([k_s], k_s+1)\cap \Gamma (S, J_{s-1}).
\end{equation}
Clearly,
\begin{equation}\label{product-al}
   \big\vert   \Gamma (S, J) \cap \D(I,H) \big\vert  =
   |\D(I,H)|\,  \prod _{s=1}^{\ell}
   \frac{|\Gamma (S, J_s)\cap \D(I,H)|}{|\Gamma (S, J_{s-1})\cap \D(I,H)|}.
\end{equation}
For  $1\leq s\leq \ell$ define
$$
  {\cal{F}}_s = \left\{F \subset [n]\times[n] \, : \,  \D([k_{s}], F)\cap \D(I,H)
   \subset \Gamma (S, J_{s-1}) \right\}.
$$
Then by  \eqref{eq-coro-split} we have
$$
     \Gamma (S, J_{s})\cap \D(I,H)= \bigsqcup_{F\in\mathcal{F}_s}
     \Gamma ([k_s], k_s+1)\cap  \D([k_{s}], F)\cap \D(I,H).
$$
Applying Lemma~\ref{lem2-graph-sparse} we obtain
\begin{align*}
  \big\vert     \Gamma (S, J_s)\cap \D(I,H) \big\vert &=  \sum _{F \in {\cal{F}}_s} | \Gamma ([k_s], k_s+1)\cap  \D([k_{s}], F)\cap \D(I,H)| \\
   &\leq \gamma_{k_s} \sum _{F \in {\cal{F}}_s} | \D([k_{s}], F)\cap \D(I,H)|\\
    &\leq    \gamma_k  | \Gamma (S, J_{s-1}) \cap \D(I,H)| ,
\end{align*}
where the last inequality follows from the definition of ${\cal{F}}_s$ and $k_s\leq k$.
This and \eqref{product-al} imply the result.
\qed


We are now ready to prove Theorem~\ref{th-graph-sparse}.
In the proof, we will use Corollary~\ref{cor-al}, together with an easy observation
that the condition $|\innbr(S)| \leq (1-\varepsilon ) d|S|$ for an (ordered) subset $S$ of vertices implies
that proportionally many  vertices in $S$ have at least $\varepsilon d/2$ common in-neighbors
with the union of the preceeding vertices.

\medskip

\noindent {\bf Proof of Theorem~\ref{th-graph-sparse}.}
Let $\nrangr\in \Gamma _k$ and $S$ be as in the definition of $\Gamma _k$. For $j\in S$ consider
$$
  A_j = \bigcup_{i\in S, i<j} \inco(i, j)
$$
and denote by $m_j$ its cardinality. Note that for
$j_0=\min \{j \, : \, j\in S\}$ one has $A_{j_0}=\emptyset$ and $m_{j_0}=0$. Note also
$$
  \lvert\innbr(S)\rvert = \sum_{j\in S} (d- m_j).
$$
Let $\delta =\eps/2$ and consider $J' :=\{ j\in S \, : \, m_j \geq \delta d\}.$  Since $m_{j_0}=0$,
$$
  (1-\eps )\, d\, |S| \geq  \lvert\innbr(S)\rvert \geq \sum_{j\in S\setminus J'} (d- m_j) > (1-\eps/2)\, d \, (|S|-|J'|),
$$
which implies
$$
   |J'| > \frac{\eps}{2-\eps}\, |S| > \frac{\eps}{2} \, |S|.
$$
Hence, for every $\nrangr\in \Gamma _k$ there exists $S\subset I^c$ with $|S|=k$,
and $J \subset S$ such that
\begin{equation*}
  |J |=\lceil \eps k /2 \rceil:=\ell \quad \mbox{ and } \quad m_j\geq \delta d \, \,
\mbox{ for all } \, \,  j\in J.
\end{equation*}
Thus
$$
  \Gamma _k \subset \bigcup _{|S|=k} \, \, \bigcup_{J\subset S, |J|=\ell}\,\,  \Gamma(S, J).
$$
By Corollary~\ref{cor-al} we have
$$
  \lvert \Gamma _k \cap \D(I,F) \rvert  \leq {n\choose k} \, {k \choose \ell} \,
  {\gamma_k}^{\ell} \, |\D(I,F) | \leq  \left(\frac{en}{k}\right)^k\,
  \left(\frac{e k}{\ell}\right)^{\ell} \, {\gamma_k}^{\ell} \, |\D(I,F)|.
$$
We assume that $\eps k \geq 2$ (the case $\eps k <2$, in which  $\ell =1$, is treated similarly). Using $\eps \geq \max\{\sqrt{\ln d /d}, \sqrt{8/d}\}$, by direct calculations we observe
$$
  \left(\frac{en}{k}\right)^k\, \left(\frac{e k}{\ell}\right)^{\ell} \,
 \left(\frac{4 e kd}{\eps n}\right)^{\eps d\ell/2}    \leq
 \left(\frac{en}{k}\right)^k\, \left(\frac{2 e}{\varepsilon}\right)^{\eps k/2} \,
 \left(\frac{4ekd}{\eps n}\right)^{\eps^2 dk /4}
  \leq \left(\frac{C_1 k d}{\eps n}\right)^{\eps^2 dk /8}
$$
%
for a sufficiently large absolute constant $C_1>0$.
Taking $c\leq 1/(3 e C_1)$, we obtain the desired estimate for $\Gamma_k$.
The second assertion of the theorem regarding $\Gamma$ follows immediately.
\qed

As we have already noted, Theorem~\ref{th-graph-sparse} essentially postulates that a
random $d$-regular digraph typically has good expansion properties in the sense that
every sufficiently small subset $S$ of its vertices has almost $d|S|$ in-neighbors and $d|S|$ out-neighbors.
In the undirected setting, expansion properties of graphs are a subject of very intense research
(see, in particular, \cite{HLW} and references therein).
As the conclusion for this subsection, we would like to recall some of the known expansion properties
of undirected random graphs and compare them with the main result of this part of our paper.

Let $G=(V,E)$ be an undirected graph on $n$ vertices. Given a subset $U \subset V$, by $\partial_V U$
we denote a set of all vertices adjacent to the set $U$ but not in $U$, i.e.
$$
 \partial_V U:=\{i\not\in U:\ \exists j\in U\, (i,j)\in E\}
  = \innbr(U)\setminus U.
$$
Similarly, let $\partial_E U$ be the set of all edges of $G$ with
exactly one endpoint in $U$.
For every $\lam\in(0,1]$, we define {\it the $\lam$-vertex isoperimetric number}
$$
 i_{\lam,V}(G):=\min_{\vert U\vert \leq \lam n} \frac{ \vert \partial_V U\vert}{\vert U\vert},
$$
and, for every $\lam\in(0,1/2]$, {\it the $\lam$-edge isoperimetric number}
$$
i_{\lam,E}(G):=\min_{\vert U\vert \leq \lam n} \frac{ \vert \partial_E U\vert}{\vert U\vert}.
$$
For $\lam=1/2$, the above quantities are simply called the vertex and the edge
isoperimetric numbers, denoted  by $i_V(G)$ and $i_E(G)$.
Since $\vert\partial_V U\vert \leq \vert \partial_E U\vert\leq d\vert\partial_V U\vert$,  for every $\lam\in(0,1/2]$ we have
\begin{equation}\label{relation-isop-edge-vertex}
i_{\lam,V}(G)\leq i_{\lam,E}(G)\leq d i_{\lam,V}(G).
\end{equation}

Now, let $\rangr$ be a $d$-regular graph uniformly distributed on the set $\UndirIntro$.
In \cite{MR947025} it was shown that for large enough fixed $d$
\begin{equation}\label{eq-bollobas}
i_E(\rangr)\geq d/2 -\sqrt{d\ln 2},
\end{equation}
with probability going to one with $n$.
This result was generalized in
\cite{MR3183044}, where is was shown that
$$
 i_{\lam,E}(\rangr)\geq d(1-\lam +o(1))
 $$
with probability going to one with $n$, where $o(1)$ depends on $d$ and
can be made arbitrarily small by increasing $d$.
Note that the relation \eqref{relation-isop-edge-vertex} together with
results from \cite{MR947025, MR3183044} immediately implies
$$i_{\lam,V}(\rangr)\geq 1-\lam +o(1)$$
(where the bound should be interpreted in the same way as before), however
the bound is far from being optimal.
An estimate for the second eigenvalue of $\rangr$ proved in \cite{MR2437174} implies that for a fixed $d$
with large probability (going to one with $n$)
$$
i_V(\rangr)\geq 1-8/d + O(1/d^2).
$$
Moreover, for every $d$ and $\delta>0$ for small enough $\lam=\lam(d, \delta)>0$ the parameter $i_{\lam,V}$
(corresponding to expansions of small subsets of $V$) can be estimated as
$$i_{\lam,V}(\rangr)\geq d-2-\delta$$
(see \cite[Theorem~4.16]{HLW}).

Our main result of this subsection can be interpreted as an expansion
property of regular digraphs for small vertex subsets. We define
the vertex isoperimetric number $i_{\lam,V}$ for digraphs by the same
formula as for undirected graphs. Theorem~\ref{th-graph-sparse} has the
following consequence, which, in particular, provides quantitative
estimates of $i_{\lam,V}$ for  $d$ growing together with $n$ to infinity.

\begin{cor}
Let $8\leq d\leq n$ and $\varepsilon\in (0,1)$. Assume that
$$
\varepsilon^2 \geq \frac{\max\{8, \ln d\}}{d}, \quad d\leq \frac{c\eps n}{2}
\quad\text{and}\quad  \lam(\varepsilon) := \frac{c\varepsilon}{d}.
$$
Further, let $\rangr$ be uniformly distributed on $\D$.
Then
$$
i_{\lam(\varepsilon), V}(\rangr) \geq (1-\varepsilon)d-1
$$
with probability at least
$$
 1- \exp\left(-\frac{\eps^2 d}{8} \ln\left(\frac{ec\eps n}{d}\right)\right).
$$
\end{cor}


\subsection{On existence of edges connecting large vertex subsets}
\label{section-independence-graph}


In this part, we consider the following problem.
Let $\rangr$ be uniformly distributed on $\D$ and let $I$ and $J$ be two (large enough) subsets of $[n]$. We want to estimate the probability that $\rangr$ has no edges connecting a vertex from $I$ to a vertex from $J$. The main result of the subsection is the following theorem.

\begin{theorem}\label{th-graph-spread}
There exist absolute  constants $c>0$ and
$C, C_1\geq 1$ such that the following holds.
Let $C_1\leq d\leq  n/24$ and let natural
numbers $\ell$ and $r$ satisfy
\begin{equation*}
 \frac{n}{4}\geq r\geq\ell \geq  \frac{C n \ln (en/r)}{d}.
\end{equation*}
 Then
$$
 \P\left\{
 \bigcup
 \DO(I,J) \right\}
 \leq \exp\left(-\frac{c r\ell d}{n}\right),
$$
where the union is taken over all $I, J\subset [n]$ with
$\vert I\vert \geq \ell$ and  $\vert J\vert \geq r$.
\end{theorem}
\begin{Remark}
Obviously, the roles of $\ell$ and $r$ in this theorem are
interchangeable and the assumptions on $\ell$ and $r$ imply
that $\ell \geq Cn/d$ and
$r\geq C_1 n \ln d/d$.
\end{Remark}
\begin{Remark}
We would like to notice that adding an assumption
$\ell \geq 4d^2$ in this theorem,
we could simplify its proof (we would not need quite technical Lemma~\ref{lem1-graph-spread} below)
\end{Remark}
\begin{Remark}
The statement of the theorem can be related to known results on the independence number of
random undirected graphs. Recall that the independence number $\alpha(\nrangr)$ of a
graph $\nrangr$ is the cardinality of the largest subset of its vertices such that no two vertices of the subset are adjacent.
Suppose now that $\rangr$ is uniformly distributed on $\UndirIntro$.
For $d\to\infty$ with $d\leq n^{\theta}$ for some fixed $\theta<1$, it was shown in \cite{MR1142268} and  \cite{CFRR}
that, as $n$ goes to infinity, the ratio $\alpha(\rangr)/\bigl(2n d^{-1}\ln d\bigr)$ converges to $1$ in probability.
Moreover, in \cite{MR1839497} it was verified that in the range $n^{\theta}\leq d\leq 0.9n$ (for a sufficiently
large $\theta<1$), the asymptotic value of $\alpha(\rangr)$ is $2\ln d/\ln (n/(n-d))$, which is equivalent to $2n\ln d/d$ when
$d/n$ is small. Taking $I=J$ in Theorem~\ref{th-graph-spread},  we observe a bound of the same order for random digraphs,
which can be interpreted as a large deviation estimate for the independence
number as follows.

\begin{cor}
 There exist absolute positive constants $c,C$ such that for every
 $2\leq d\leq n/24$ and a random digraph $\rangr$ uniformly distributed
 on $\D$ one has
\begin{equation*}
 \P\left\{\alpha(\rangr)>  C\, \frac{n\ln d}{d}\right\} \leq
 \exp\left(-\frac{c n\ln^2d}{d}\right).
\end{equation*}
\end{cor}
%
\end{Remark}

\medskip

We first give an idea of the proof of Theorem~\ref{th-graph-spread}. Fix
two sets of vertices $I$ and $J$ of sizes $\ell$ and $r$. Our strategy is
to start with two small subsets of $I$, $J$ and to arrive to
$I$, $J$ by adding one vertex at a time.
Suppose that $I_1\subset I$ and $J_1\subset J$
and $S$ is a subset of graphs from $\D$ with no edges departing from $I_1$ and landing in $J_1$. We
add a vertex from $J\setminus J_1$ to $J_1$ to form a set $J_2$ and check whether
the property of having no edges connecting $I_1$ to $J_2$ is preserved, using the simple switching.
More precisely, when conditioning on the set of graphs $S$,
we estimate the proportion of graphs in $S$ such that there are no edges departing from $I_1$
to the vertex added. We perform an analogous procedure by adding a vertex to $I_1$ and
continue until the whole sets $I$ and $J$ are reconstructed.

Note that a similar argument can be applied in the undirected setting to estimate probability of large
deviation for the independence number (the sets $I$ and $J$ shall be equal in this situation).
We omit the proof of the undirected case as it is a simple adaptation of the argument
of Theorem~\ref{th-graph-spread} and is not of interest in the present paper.

We start with a lemma which can be described as follows:
given two sets of vertices $[p]$ and $[k]$, among
graphs having no edges departing from $[p]$ to $[k]$, we count
how many have no departing edges from $[p]$ to the vertex $k+1$.
The proof of Theorem~\ref{th-graph-spread} will then follow by iterating this lemma.

\begin{lemma}\label{lem1-square-spread-graph}
Let $20\leq d\leq n/24$ and  $4e^2n/d \leq  p, k \leq  n/4$. Then
$$
 \max\left\{ \vert \DO([p],[k+1])\vert,\, \, \vert \DO([p +1],[k])\vert \right\} \leq   \exp\left(-\frac{pd}{4e^2n}\right)\, \vert \DO([p],[k])\vert.
$$
\end{lemma}

To prove this lemma we need the following rather technical statement,
which shows that for most graphs under consideration every
two vertices have a relatively small number of common out-neighbors.
For reader's convenience we postpone its proof to the end of this section.

\begin{lemma}\label{lem1-graph-spread}
Let $\varepsilon\in (0,1)$, $0\leq k\leq n/4$, $0\leq p\leq n$, and  $d\leq \varepsilon n/12$.
Then
$$
 \lvert \DO([p], [k])\setminus \Dco(\varepsilon) \rvert
 \leq \frac{n^2}{2}\, \left(\frac{2e d}{\varepsilon n}\right)^{\varepsilon d}\,
 \lvert \DO ([p], [k])\rvert ,
$$
where $\DO([p], [k])=\D$ if $p=0$ or $k=0$.
\end{lemma}

\medskip

\noindent {\bf Proof of Lemma~\ref{lem1-square-spread-graph}.}
We prove the bound for $\DO([p],[k+1])$, the other bound is obtained by passing to the transpose graph.

Fix  $q:= \lceil pd/(2e^2n)\rceil $. Denote
$$
 Q := \DO([p],[k+1])\cap \Dco(1/2)
$$
and 
$$
 Q(q):= \{\nrangr\in \DO([p],[k]):\ \vert \edg([p], k+1)\vert =q\}.
$$

To estimate cardinalities we construct a relation $R$ between $Q$ and $Q(q)$. We say that $(\nrangr,\nrangr')\in R$ for some 
$\nrangr\in Q$ and $\nrangr' \in Q(q)$ if
$\nrangr'$ can be obtained
from $\nrangr$ using the simple switchings as follows. First
choose $q$ vertices $1\leq v_1<v_2<\ldots<v_q\leq p$. There are 
$p \choose q$ such choices. Then choose $q$ edges in  
$\edg([p]^c, k+1)$, say $(w_i, k+1)$, $i\leq q$, with $p<w_1<w_2<\ldots<w_q\leq n$.  There are 
$d \choose q$ such choices. Finally for every $i\leq q$ choose 
$$
 j(i)\in\outnbr(v_i)\setminus \outnbr(w_i).
$$ 
Since $\nrangr\in \Dco(1/2)$, for every $i\leq q$ there are at least $d/2$ choices
of $j(i)$. For every $i\leq q$ we destroy edges $(w_i, k+1)$,
$(v_i, j(i))$ and create edges $(v_i, k+1)$,
$(w_i, j(i))$. We have  
\begin{equation}\label{eq-image-lem-th3}
 \lvert R (\nrangr)\rvert \geq {p\choose q}\, {d\choose q}
  \, \left(\frac{d}{2}\right)^q \geq 
   \left(\frac{p d^2}{2q^2}\right)^q. 
\end{equation}

Now we estimate the cardinalities of preimages.
Let $\nrangr' \in R(Q)$. We bound $|R^{-1}(\nrangr')|$ from above.
To reconstruct a possible $\nrangr\in  Q$ with  $(\nrangr, \nrangr')\in R$, we perform simple switchings as follows.
Write $\edgpr([p],k+1)$ as $(v_1, k+1), \ldots, (v_q, k+1)$ with 
$1\leq v_1<\ldots < v_q\leq p$.  Choose $q$ vertices 
$p<w_1<\ldots w_q\leq n$ such that 
$$
 w_i\in [p]^c\setminus \innbrpr(k+1)
$$
for all $i\leq q$. There are 
$$
 {n-p-(d-q) \choose q} \leq \left(\frac{e n}{q}\right)^q
$$ 
such choices. For every $i\leq q$ find 
$$
 j\in \left(\outnbrpr(w_i) \cap [k+1]^c \right)\setminus
 \outnbrpr(v_i) 
$$
(there are at most $d$ such choices). 
For every $i\leq q$ we destroy edges $(v_i, k+1)$,
$(w_i, j(i))$ and create edges $(w_i, k+1)$,
$(v_i, j(i))$. We obtain 
\begin{equation*}
 \lvert R^{-1}(\nrangr')\rvert \leq 
 \left(\frac{e n}{q}\right)^q \, d^q.
 \end{equation*} 
Claim~\ref{multi-al} together with the last bound and
\eqref{eq-image-lem-th3} yields
$$
\lvert Q\rvert
\leq \left(\frac{2 e n q}{pd}\right)^q
\lvert  Q(q)\rvert.
$$
By the choice of $q$, we have 
$\lvert Q\rvert \leq \exp(-pd/(2e^2n)) \lvert  Q(q)\rvert.$
This, together with Lemma~\ref{lem1-graph-spread}, implies  
\begin{align*}
\vert\DO([p],[k])
\vert
&\geq \vert Q(q)\vert \geq \exp\left(\frac{pd}{2 e^2 n}\right)\,
\vert  \DO([p],[k+1])\cap \Dco(1/2) \vert \\
&= \exp\left(\frac{pd}{2 e^2 n}\right)\,
\left( \vert  \DO([p],[k+1]) \vert  -  
\vert  \DO([p],[k+1])\setminus \Dco(1/2) \vert 
 \right)\\
&
\geq \exp\left(\frac{pd}{2 e^2 n}\right)\, 
\left(1-\frac{n^2}{2}\, \left(\frac{4 e d}{n}
\right)^{d/2}\right) \,
\vert  \DO([p],[k+1]) \vert ,
\end{align*}
which implies the desired result.
\qed


\medskip \noindent {\bf Proof of Theorem~\ref{th-graph-spread}.}
It is enough to prove the theorem for the union over all $|I|=\ell$ and $|J|=r$.
By the union bound, we have
\begin{align}\label{eq-start-graph-spread}
  \P\left\{  {\bigcup}\, \DO(I,J) \right\}
  \leq {n\choose \ell} \, {n\choose r} \,
    \frac{\lvert\DO([\ell], [r])\rvert}{\lvert\D\rvert}
    \leq \left(\frac{en}{r}\right)^{2r}
    \,  \frac{\lvert\DO([\ell], [r])\rvert}{\lvert\D\rvert}.
\end{align}
Setting $\DO([0], [0])=\D$ and using $\DO([k], [k])\supset \DO([k+1], [k+1])$, we get
\begin{align}
    \frac{\lvert\DO([\ell], [r])\rvert}{\lvert\D\rvert}
   &= \prod_{k=0}^{\ell-1}\,  \frac{ \lvert\DO([k+1],
    [k+1])\rvert}{\lvert\DO([k], [k])\rvert}
    \,\, \prod_{k=\ell}^{r-1}\,  \frac{ \lvert\DO([\ell],
    [k+1])\rvert}{\lvert\DO([\ell], [k])\rvert}\nonumber\\
   &\leq \prod_{k=\lceil\ell/2\rceil}^{\ell-1}\,  \frac{ \lvert\DO([k+1],
    [k+1])\rvert}{\lvert\DO([k], [k])\rvert}
    \,\, \prod_{k=\ell}^{r-1}\,  \frac{ \lvert\DO([\ell],
    [k+1])\rvert}{\lvert\DO([\ell], [k])\rvert}.\label{eq-th3-start-product}
\end{align}
Further, we write
$$
\frac{ \lvert\DO([k+1],
    [k+1])\rvert}{\lvert\DO([k], [k])\rvert} =\frac{ \lvert\DO([k+1],
    [k+1])\rvert}{\lvert\DO([k], [k+1])\rvert} \cdot \frac{ \lvert\DO([k],
    [k+1])\rvert}{\lvert\DO([k], [k])\rvert},
$$
and applying Lemma~\ref{lem1-square-spread-graph},
for every $\lceil \ell/2\rceil \leq k\leq \ell-1$ we observe
\begin{equation*}\label{eq-th3-product-part1}
\frac{ \lvert\DO([k+1],
    [k+1])\rvert}{\lvert\DO([k], [k])\rvert}
  \leq \exp\left(-\frac{kd}{2 e^2 n}\right),
  \end{equation*}
and for every $\ell\leq k\leq r-1$,
\begin{equation*}\label{eq-th3-product-part2}
\frac{ \lvert\DO([\ell],
    [k+1])\rvert}{\lvert\DO([\ell], [k])\rvert}
  \leq \exp\left(-\frac{\ell d}{4 e^2 n}\right).
  \end{equation*}
Thus \eqref{eq-th3-start-product} implies
\begin{equation*}
 \frac{\lvert\DO([\ell], [r])\rvert}{\lvert\D\rvert}
 \leq \exp\left(-\frac{\ell rd}{8 e^2 n}\right).
\end{equation*}
Combining this bound and (\ref{eq-start-graph-spread}) and using
that $\ell \geq  C n \ln (en/r) /d$ we complete the proof.
\qed

\medskip

\noindent {\bf Proof of Lemma~\ref{lem1-graph-spread}.}
Clearly,
\begin{equation*}
  \lvert \DO([p], [k])\setminus \Dco(\varepsilon)\rvert
  \leq \sum_{i< j} \lvert \DO([p], [k])\setminus \D_{i,j}^{co}(\varepsilon)\rvert .
\end{equation*}

Fix $1\leq i<j\leq n$. For $0\leq q\leq d$, denote
$$
  Q(q):=\bigl\{ \nrangr \in \DO([p], [k]):\, \vert \outco(i,j)\vert =q\bigr\}
$$
and
$$
  Q:=\DO([p], [k])\setminus \D_{i,j}^{co}(\varepsilon)=
  \bigsqcup _{q=\lfloor\varepsilon d\rfloor + 1}^d Q(q).
$$

First, for every $1\leq q\leq d$ we will compare the cardinalities of $Q(q)$ and $Q(q-1)$.
To this end, we will define a relation $R_q$ between the sets $Q(q)$ and $Q(q-1)$ in the following way.

Let $\nrangr \in Q(q)$. Then there exist
$j_1<...<j_{q}$ such that for every $\ell\leq q$
$$
j_{\ell}\in \outnbr(i)\cap \outnbr(j).
$$
Note that for every $\ell \leq q$,  there are $d^2$ edges inside
$
\outedg\big(\innbr(j_\ell)\big).
$
Also, there are at least $(n-k-(2d-q))d$ edges in $\inedg\big([k]^c\setminus \outnbr(\{i,j\})\big)$.
Therefore, for  $\ell \leq q$, the set
$$
  E_\ell:=\edg\big([n]\setminus \innbr(j_\ell), [k]^c\setminus \outnbr(\{i,j\})\big)
$$
is of cardinality at least
$$
  |E_\ell|\geq (n-k-(2d-q))d-d^2\geq   nd/2.
$$

We say that $(\nrangr, \nrangr') \in R_q$ for some $\nrangr' \in Q(q-1)$
if $\nrangr'$ can be obtained from $\nrangr$ in the following way. First
we choose $\ell \leq q$ and an edge $(u,v)\in E_\ell$. Note $v\in [k]^c$
and $u\ne i$. Since $v\not\in \outnbr(j)$ then we can destroy the edge
$(j,j_\ell)$ and create the edge $(j,v)$. Since $u\not\in\innbr(j_\ell)$
then we can destroy the edge $(u,v)$ and create the edge $(u,j_\ell)$.
Thus, we obtain $\nrangr'$ by the simple switching on vertices $u,v,j,j_\ell$.
It is not difficult to see that we have not created any edges between $[p]$ and $[k]$,
hence $\nrangr'$ indeed belongs to $Q(q-1)$.
Counting the admissible simple switchings, we get for every $\nrangr\in Q(q)$,
\begin{equation}\label{eq-image-Rq-lem5}
\vert R_q(\nrangr)\vert\geq \frac{qnd}{2}.
\end{equation}

Now we estimate the cardinalities of  preimages.
Let $\nrangr' \in R_q\bigl( Q(q)\bigr)$.
In order to reconstruct a possible $\nrangr$ for which
$(\nrangr, \nrangr')\in R_q$, we need to perform the simple switching which removes an edge $(j, v)$
with $v\not\in  \outnbrr(i)$ and recreates an edge $(j, w)$ for some
$$
 w\in  \outnbrr(i)\setminus \outnbrr(j).
$$
There are at most $d-q+1$ choices for such $v$ and
at most $d-q+1$ choices  for such $w$.
For the second part of the switching, we have
at most $d$ possible choices. Therefore,
\begin{equation*}
 \lvert R_q^{-1}(\nrangr')\rvert \leq d(d-q+1)^2\leq d^3.
\end{equation*}

Using this bound, \eqref{eq-image-Rq-lem5}, and Claim~\ref{multi-al},
we obtain that
$$
\lvert Q(q)\rvert \leq \left(\frac{2 d^2}{qn}\right)\cdot  \lvert Q(q-1)\rvert
$$
and, applying this successively,
\begin{equation*}
    \lvert Q(q)\rvert     \leq          \left(\frac{2 d^2}{n}\right)^{q} \, \,
    \frac{1}{q!}\, \,   \lvert Q(0)\rvert.
\end{equation*}
Since
         $q!\geq 2(q/e)^q$ and $2e d/(\varepsilon n)\leq 1/2$, this implies
$$
  \vert Q\vert= \sum_{q=\lfloor\varepsilon d\rfloor + 1}^d\lvert Q(q)\rvert
  \leq \frac{1}{2}\,\,  \sum_{q=\lfloor\varepsilon d\rfloor + 1}^d
  \left(\frac{2e d}{\varepsilon n}\right)^{q} \, \lvert Q(0)\rvert
  \leq \left(\frac{2e d}{\varepsilon n}\right)^{\varepsilon d}\,
  \lvert Q(0)\rvert .
$$
Using that $Q(0)\subset \DO([p],[k])$ and that there are $n(n-1)/2$ pairs $i<j$, we obtain the desired result.
\qed

\subsection{An anti-concentration property for random digraphs}
\label{section-anticoncentration-graph}

For every $\nrangr\in \D$, $J\subset [n]$ and $i\in[n]$,
we define $\delta_i^{J}\in \{0,1\}$ by
$$
\delta_i^J=\delta_i^J(\nrangr) :=\left\{
	\begin{array}{ll}
		1  & \mbox{if } i\in \innbr(J), \\
		0 & \mbox{otherwise}.
	\end{array}
\right.
$$
Denote $\delta^J:=(\delta_1^J,\ldots, \delta_n^J)\in \{0,1\}^n$.
The vector $\delta^J$ can be regarded as an indicator of the vertices that are connected to $J$,
without specifying how many edges connect a vertex to $J$.

 Taking $v\in\{0,1\}^n$ and conditioning on the realization $\delta^J=v$,
we obtain a class of graphs with a particular arrangement of the edges. Namely,
if a vertex $i$ of a graph in the class is not connected to $J$ then all graphs in the class have the same property.  In this section we estimate the cardinalities of such classes  generated by vertices of the cube, under additional assumption that a part of the graph is ``frozen." We show that if the size of the set $J$ is at most $cn/d$ then a large proportion of such classes are ``approximately'' of the same size. In other words, we prove that  the distribution of $\delta^J$ is similar to that of a random vector uniformly distributed on the discrete cube $\{0,1\}^n$ in the sense that for each fixed $v\in\{0,1\}^n$ the probability that $\delta^J=v$ is exponentially small. This makes a link to the anti-concentration results in the Littlewood-Offord theory. We start with a simplified version of this result, when there is no ``frozen" part. In this case it is a rather straightforward consequence of Theorem~\ref{th-graph-sparse}.


\begin{prop}\label{claim-anti}
Let $8\leq d\leq n$ and $J\subset[n]$. Let $v\in\{0,1\}^n$ and $m:=|\supp v|$.
Then
\begin{equation*}\label{P-ev-I1}
\P\{ \delta^J=v\}\leq \binom{n}{m}^{-1}\le
  \exp\left(-m\ln\frac{ n}{m}\right).
\end{equation*}
Moreover, if $|J|\leq c n/d$, then
\begin{equation*}
\label{P-ev-I0}
\P\{ \delta^J=v\}\leq
  \exp\left(-cd|J|\ln\frac{ cn}{d|J|}\right),
\end{equation*}
where $c$ is an absolute positive constant.
\end{prop}

\begin{Remark}
Note that  $\max \{d, |J|\}\leq|\supp \delta^J|\le d|J|$, therefore $\P\{ \delta^J=v\}=0$ unless $\max \{d,|J|\}\leq m \leq d|J|$.
\end{Remark}

\proof
 Without loss of generality we assume  that $\max\{d,|J|\}\leq m \leq d|J|$.
 Consider the following subset of the discrete cube
$$
T=\{w\in\{0,1\}^n \, :\, |\supp w|=m\}.
$$
Clearly, every $w\in T$ can be obtained by a permutation of
the coordinates of $v$.
Since the distribution of a random graph is invariant under permutations,
$\P\{ \delta^J=v\}=\P\{ \delta^J=w\}$ for every $w\in T$. Therefore,
$$
\p\{\delta^J=v\}\le |T|^{-1}=\binom{n}{m}^{-1}\le  \exp\left(-m\ln\frac{n}{m}\right),
$$
which proves the first bound and the ``moreover" part in the case
$m\geq d|J|/2$.

Suppose now that $\vert J\vert \leq c n/d$ and $m\leq d|J|/2$.
Applying Theorem~\ref{th-graph-sparse} with $S=J$, $I=\emptyset$, and $\varepsilon=1/2$,
we observe
$$
\P\{ |\supp \delta^J|\le d|J|/2 \}\leq
  \exp\left(-cd|J|\ln\frac{ cn}{d|J|}\right),
$$
which completes the proof of the ``moreover" part.
\qed

Now we turn to the main theorem of this section, which will play a key role in the ``matrix" part of our paper.  We obtain an anti-concentration property for the vector $\delta^J$  even under an  assumption that a part of edges is  ``frozen.'' It requires a more delicate argument.


\begin{theorem}\label{th-anticoncentration}
There exist two absolute positive constants $c, \tilde c$ such that
the following holds.
Let $32\leq d\leq cn$ and let $I, J$ be disjoint subsets of $[n]$ such that
\begin{equation*}\label{eq-choiceI-th-anticoncentration}
 \vert I\vert \leq \frac{d\vert J\vert}{32} \quad \text{ \ and \ }\quad
 8\leq \vert J \vert  \leq \frac{8 c n}{d} .
\end{equation*}
Let $F\subset [n]\times[n]$ be such that $\D(I,F)\neq\emptyset$ and let $v\in \{0,1\}^n$. Then
$$
  \P\{ \delta^J =v  \mid  \D(I,F)\}\leq
  2\exp\left(- \tilde c d\vert J\vert \ln\left(\frac{n}{d\vert J\vert}\right)\right).
$$
\end{theorem}



To prove this theorem we first estimate the size
of the class of graphs given by a realization of a subset of coordinates of $\delta^J$.
More precisely, restricted to a subset of graphs with predefined  out-edges for
the first $i-1$ vertices of $\delta^J$, we count
the number of graphs for which the vertex $i$ is connected to $J$.
 In other words, conditioning on the realization of the first
$i-1$ coordinates of $\delta^J$, we estimate the probability that $\delta_i^J=1$.
In Lemma \ref{lem-anticoncentration} below we show that this probability is of the order $d\vert J\vert /n$ for a wide range of $i$-s. In a sense, this shows that
the sets of out-edges restricted on $J$ for different vertices  behave like independent. Indeed, in the Erd\H{o}s--R\'enyi model, when the edges are distributed independently
with probability of having an edge equals $d/n$, the probability that a vertex $i$ is connected to $J$ is of order $d\vert J\vert /n$.


We need the following notation. For $\eps \in (0,1)$ and $J\subset [n]$ let
$$
 \Lambda(\varepsilon ,J)= \big\{\nrangr\in \D:\ \vert \innbr(J)\vert
 \geq (1-\varepsilon) d\vert J\vert \big\} .
$$

\medskip

\begin{lemma}\label{lem-anticoncentration}
Let $2\leq d\leq n/32$.
Let $F\subset[n]\times[n]$ and $I, J$ be disjoint subsets of $[n]$ satisfying
\begin{equation}\label{eq-choiceI-lem-anticoncentration}
 \vert I\vert \leq  \frac{d\vert J\vert}{32} \quad \text{ \ and \ }
 \quad 8 \leq \vert J \vert  \leq \frac{n}{4d}.
\end{equation}
Then there exists a permutation $\sigma\in \Pi_n$ such that for every
$$
  2\vert I\vert\leq i_1< i_2< \ldots < i_{d\vert J\vert/16},
$$
every $s \leq d\vert J\vert/16$ and $\mathcal{H}\subset 2^{[n]\times[n]}$
satisfying
$$
  {\widetilde{\Gamma}}:= \left\{ \nrangr\in \D(I,F)\cap \Lambda \left(\tfrac{1}{2}, J\right)\, : \,
  \edg\big(\sigma([2\vert I\vert]\cup\{i_1,\ldots,i_{s-1}\}), I^c\big)\in \mathcal{H} \right\} \ne \emptyset
$$
one has
\begin{equation*}\label{eq-lem-anticoncentration}
 \frac{d\vert J\vert}{9n}\leq \P \left\{ \delta_{\sigma(i_s)}^J =1   \mid {\widetilde{\Gamma}}  \right\}
 \leq  \frac{2d\vert J\vert}{n}.
\end{equation*}
\end{lemma}

\medskip

As the proof of lemma is rather technical, we postpone it to the end
of this section.

\medskip
\noindent{\bf Proof of Theorem~\ref{th-anticoncentration}. }
Fix  $F\subset [n]\times[n]$ with $\D(I,F)\neq\emptyset$ and $v\in \{0,1\}^n$.
Let $\sigma$ be a permutation given by Lemma~\ref{lem-anticoncentration}.

Denote $B:= \D(I,F)\cap \Lambda \left(\tfrac{1}{2}, J\right)$.
Since $J\subset I^c$, applying Theorem~\ref{th-graph-sparse} with $\eps =1/2$ and $k=|J|$,
we get that for some appropriate constant $\tilde c$
\begin{equation*}\label{new-eq-13}
 \P\left\{  \Lambda \left(\tfrac{1}{2}, J\right) \mid \D(I,F) \right\}
 \geq 1-\exp\left(- \tilde c d\vert J\vert \ln\Big(\frac{n}{d\vert J\vert}\Big)\right),
\end{equation*}
which in particular implies that  $B$ is non-empty.

Using this we have
\begin{align*}
 \P\{ \delta^J =v  \mid  \D(I,F) \} &\leq \P\left\{ \delta^J =v  \mid
 \D(I,F)\cap \Lambda \left(\tfrac{1}{2}, J\right) \right\} + \P\left\{  \Lambda^c
\left(\tfrac{1}{2}, J\right)\mid  \D(I,F) \right\}\nonumber\\
 & \leq \P\left\{ \delta^J =v  \mid
 B\right\}+ \exp\left(- \tilde c d\vert J\vert \ln\left(\frac{n}{d\vert J\vert}\right)\right).\label{eq-anti-new-th}
\end{align*}

Therefore, it is enough to estimate the first term in the previous inequality.
Note that if $\vert\supp v\vert < d\vert J\vert/2$ then
$$
\P\left\{ \delta^J =v  \mid
 \D(I,F)\cap \Lambda \left(\tfrac{1}{2}, J\right) \right\}=0.
$$
Assume that $\vert\supp v\vert\geq d\vert J\vert/2$ and denote $m=d|J|/16$.  Since
$2\vert I\vert\leq  m$, there exist
$$
 2\vert I\vert \leq i_1< i_2< \ldots< i_{m}
$$
such that for every $s\leq  m$ one has $v_{\sigma(i_s)}=1.$
Let $Q_1=[2\vert I\vert]$. For every $2\leq s\leq m+1$, define
$Q_s:=Q_1\cup\{i_1,\ldots,i_{s-1}\}$ and
$$
  \mathcal{H}_s:=\bigl\{H\subset \sigma(Q_s)\times I^c
  \, :\,
   \forall \ell \in Q_s\, \, \, \, \,
   ``v_{\sigma(\ell)}=0"
   \Leftrightarrow ``\forall j\in J \, \, \, \, (\sigma(\ell),j)\not\in H"\bigr\}.
$$
In words, $\mathcal{H}_s$ is the collection of all possible realizations of configurations of edges connecting $\sigma(Q_s)$
to $I^c$, such that $\sigma(\ell)$ is not connected to $J\subset I^c$ if and only if $v_{\sigma(\ell)}=0$ ($\ell\in Q_s$).
Note that
\begin{equation*}
  A_s:= \big\{\nrangr \in \D\, : \, \forall \ell \in Q_s\quad \delta^J_{\sigma(\ell)}=v_{\sigma(\ell)} \big\} =
  \big\{\nrangr \in \D \, : \,  \edg\big(\sigma(Q_s), I^c\big)\in \mathcal{H}_s\big\}.
\end{equation*}
Denote
$$
  B_s:= \big\{\nrangr \in \D\, : \,
  \delta^J_{\sigma(i_s)}=1 \big\}.
$$
Since $v_{\sigma(i_s)}=1$ for every $s\leq m$ then $A_{s+1}= B_s \cap A_s$ and
$$
   \P\{ A_{s+1}  \mid B\} = \P\{  B_s \cap A_s \mid B\}
   = \P\{  B_s  \mid B\cap A_s \}\, \P\{ A_{s}  \mid B\}.
$$
Therefore,
$$
  \P\{ \delta^J =v \mid  B\}\leq  \P\{ A_{m+1}  \mid B\}
  \leq  \prod_{s = 1}^{m} \P\{  B_s  \mid B\cap A_s \} .
$$
By the assumptions of the theorem and
Lemma~\ref{lem-anticoncentration},
 for every $s\leq m$ we have
$$
  \P\{  B_s  \mid B\cap A_s \}\leq  \frac{2d|J|}{n},
$$
which implies
\begin{equation*}
\label{eq-anticoncentration-finalestimate}
  \P\{ \delta^J =v \mid  B\}\leq \left( \frac{2d|J|}{n}\right)^m\leq
  \exp\left(-\frac{d |J|}{16}\ln\left(\frac{n}{2d\vert J\vert}\right)\right) .
\end{equation*}
This completes the proof.
\qed

\medskip
It remains to prove Lemma~\ref{lem-anticoncentration}. To get the
lower bound we employ the simple switching to graphs whose $i$-th
vertex is not connected to $J$ and transform them into graphs with
the $i$-th vertex connected to $J$. To get the upper bound, we do
the opposite trick to transform graphs with only one edge relating
vertex $i$ to $J$ to a graph with no connections from $i$ to $J$. Then
we show that if $i$ is connected to $J$, it is more likely that the
number of corresponding out-edges is small. This is very natural if
we have in mind the result proven in Theorem~\ref{th-graph-sparse}.
Indeed, if vertices of a graph had a large number of out-edges
connecting them to $J$, then the number of in-neighbors to $J$
would be small, while Theorem~\ref{th-graph-sparse} states that
$\innbr(J)$ is rather large.

\medskip

\noindent {\bf Proof of Lemma~\ref{lem-anticoncentration}. }
Let $\sigma$ be a permutation such that the sequence
$$
   \big\{\vert \outnbr\big(\sigma(\ell)\big)\cap I\vert \big\}_{\ell=1}^{n}
$$
is non-increasing. Note that $\sigma$ depends only on $F$ when $\nrangr\in \D(I,F)$.

First we note that for every
 $\nrangr\in  \D(I,F)$
\begin{equation}\label{eq-many-1-out-of-I}
\forall i\geq 2\vert I\vert\, \, \, \quad  \big\vert \outnbr\big(\sigma(i)\big)\cap I^c \big\vert \geq d/2 .
\end{equation}
Indeed, otherwise there would exist $\nrangr\in  \D(I,F)$ and $i_0\geq 2\vert I\vert$ such that
$$
 \big\vert \outnbr\big(\sigma(i_0)\big)\cap I\big\vert > d/2.
$$
Since $\big\{\vert \outnbr\big(\sigma(\ell)\big)\cap I\vert\big\}_{\ell\leq n}$ is non-increasing,
then for every $\ell \leq i_0$
we would have
$$
 \big\vert \outnbr\big(\sigma(\ell)\big)\cap I \big\vert > d/2.
$$
This would imply
$$
 \big\vert \edg(\sigma([i_0]), I)\big\vert > i_0 d/2
 \geq d\vert I\vert,
$$
which is impossible.

Fix $s\leq d\vert J\vert/16$.
For  $0\leq k\leq p:=\min\{d,\vert J\vert\}$ denote
$$
 \widetilde{\Gamma}_k:= \big\{\nrangr\in \widetilde{\Gamma}:\ \big\vert \edg(\sigma(i_s),J)\big\vert =k\big\}.
$$
Clearly,
$
\widetilde{\Gamma}=\bigsqcup_{k\leq p} \widetilde{\Gamma}_k.
$

The statement of the lemma is equivalent to the following estimate
\begin{equation}\label{eq-lem-anticoncentration-goal}
 \frac{d\vert J\vert}{9n}\, \big\vert \widetilde{\Gamma}\big\vert \leq
  \big\vert \widetilde{\Gamma}\setminus \widetilde{\Gamma}_0\big\vert
 \leq \frac{2 d\vert J\vert}{n} \, \big\vert \widetilde{\Gamma}\big\vert .
\end{equation}

We first show that
\begin{equation}\label{eq-lem-anticoncentration-claim1}
 \big\vert \widetilde{\Gamma}_0\big\vert \leq \frac{8n}{d\vert J\vert}\, \big\vert \widetilde{\Gamma}_1\big\vert .
\end{equation}
Note that (\ref{eq-lem-anticoncentration-claim1}) implies the left hand side of
(\ref{eq-lem-anticoncentration-goal}).
Indeed, since  $\widetilde{\Gamma}_1\subseteq  \widetilde{\Gamma}\setminus \widetilde{\Gamma}_0$, then
(\ref{eq-lem-anticoncentration-claim1}) yields   that
$$
 \big\vert \widetilde{\Gamma}_0\big\vert \leq \frac{8n}{d\vert J\vert}\,
  \big\vert  \widetilde{\Gamma}\setminus \widetilde{\Gamma}_0\big\vert.
$$
Adding $\big\vert  \widetilde{\Gamma}\setminus \widetilde{\Gamma}_0\big\vert$ to both sides we obtain
 the left hand side of (\ref{eq-lem-anticoncentration-goal}).

In order to prove (\ref{eq-lem-anticoncentration-claim1}), we define a relation $R$
between the sets $\widetilde{\Gamma}_0$ and $\widetilde{\Gamma}_1$.
Let $\nrangr\in \widetilde{\Gamma}_0$. Since $G\in \Lambda \big(\frac{1}{2}, J\big)$
and $2\vert I\vert +s \leq d\vert J\vert/8$, then
\begin{equation}\label{eq-lem-anticoncentration-innbrJ}
\big\vert \innbr(J)\setminus \sigma([2\vert I\vert]\cup \{i_1,\ldots,i_{s-1}\})\big\vert \geq \frac{3d\vert J\vert }{8}.
\end{equation}
Denote
$$
 T:=\big(\innbr(J)\setminus \sigma([2\vert I\vert]\cup \{i_1,\ldots,i_{s-1}\})\big)\times (\outnbr\big(\sigma(i_s)\big)\cap I^c).
$$
Since $\nrangr\in \widetilde{\Gamma}_0$, that is $\big\vert \edg(\sigma(i_s),J)\big\vert =0$, we have
\begin{equation*}\label{eq-lem-anticoncentration-edgT}
\big\vert\edg\big(T\big)\big\vert \leq (d-1)\big\vert \outnbr\big(\sigma(i_s)\big)\cap I^c\big\vert .
\end{equation*}
 This together with (\ref{eq-many-1-out-of-I}), (\ref{eq-lem-anticoncentration-innbrJ}), and
$\vert J\vert \geq 8$ implies that
the set $S:=T\setminus \edg(T)$ satisfies
\begin{equation}\label{eq-lem-anticoncentration-sizeS}
  \vert S\vert \geq \left(\frac{3d\vert J\vert }{8}-d+1\right) \cdot \big\vert
  \outnbr\big(\sigma(i_s)\big)\cap I^c \big\vert
  \geq \frac{d^2\vert J\vert}{8}.
\end{equation}

  We say that $(\nrangr,\nrangr')\in  R$ for some $\nrangr'\in \widetilde{\Gamma}_1$ if $\nrangr'$ can be obtained from
$G$ in the following way.
  First choose $(v,j)\in S$. Since $j\in \outnbr\big(\sigma(i_s)\big)\cap I^c$ and
  $(v,j)\not\in \edg(T)$ then we can destroy the edge $(\sigma(i_s),j)$ and create the edge $(v,j)$.
  Since $v\in \innbr(J)$, there is $j'\in J$ such that $(v, j')$ is an edge in $\nrangr$. Since
$\nrangr\in \widetilde{\Gamma}_0$, $(\sigma(i_s), j')\not\in \nrangr$. Thus we can destroy the edge
$(v,j')$ and create the edge $(\sigma(i_s),j')$, completing the simple switching.
By \eqref{eq-lem-anticoncentration-sizeS} we get
\begin{equation*}\label{eq-imageR1-anticoncentration}
\vert  R(\nrangr)\vert\geq \frac{d^2\vert J\vert}{8}.
\end{equation*}
Note that the above transformation of $\nrangr$ does not decrease
$\vert \innbr(J)\vert$ which guarantees that $\nrangr'\in  \Lambda \big(\frac{1}{2}, J\big)$.

Now we estimate the cardinalities of preimages. Let $\nrangr'\in R(\widetilde{\Gamma}_0)$.
In order to reconstruct a possible $\nrangr$ for which $(\nrangr,\nrangr')\in  R$,
destroy the only edge $(\sigma(i_s), j')$ in $\edgpr(\sigma(i_s), J)$
and create an edge $(\ell, j')$ for $\ell \not\in \sigma([2\vert I\vert]\cup\{i_1,\ldots,i_{s-1}\})$.
There are at most $n-2\vert I\vert -(s-1)\leq n$ possible choices at this step.
To complete the simple switching, we destroy one of the edges in $\edgpr(\ell, J^c\cap I^c)$
and create an edge connecting $\sigma(i_s)$ to $J^c\cap I^c$.
There are at most $d$ possible choices here. Therefore,
\begin{equation*}\label{eq-preimageR1-anticoncentration}
\vert  R^{-1}(\nrangr')\vert\leq nd.
\end{equation*}
By Claim~\ref{multi-al}, this implies the inequality \eqref{eq-lem-anticoncentration-claim1}.

We  now show that for every $k\in\{1,\ldots, p\}$, one has
 \begin{equation}\label{eq-lem-anticoncentration-claim2}
              \vert \widetilde{\Gamma}_k\vert\leq \frac{2d\vert J\vert}{k n}\, \vert \widetilde{\Gamma}_{k-1}\vert .
 \end{equation}
Note that (\ref{eq-lem-anticoncentration-claim2}) implies the right hand side of
(\ref{eq-lem-anticoncentration-goal}).
Indeed, by (\ref{eq-lem-anticoncentration-claim2}),
$$
  \vert \widetilde{\Gamma}_k\vert\leq \left(\frac{2d\vert J\vert}{n}\right)^k \, \frac{1}{k!}\, \, \vert
  \widetilde{\Gamma}_{0}\vert,
$$
hence
$$
 \vert \widetilde{\Gamma}\vert = \vert \widetilde{\Gamma}_0\vert +\sum_{k=1}^p  \vert \widetilde{\Gamma}_k\vert
 \leq  \exp\left(\frac{2d\vert J\vert}{n}\right) \vert \widetilde{\Gamma}_0\vert ,
$$
which implies
$$
  \big\vert \widetilde{\Gamma}\setminus \widetilde{\Gamma}_0\big\vert
 \leq
 \big\vert \widetilde{\Gamma}\big\vert-\exp\left(-\frac{2 d\vert J\vert}{n}\right) \, \big\vert \widetilde{\Gamma}\big\vert
 \leq \frac{2 d\vert J\vert}{n} \, \big\vert \widetilde{\Gamma}\big\vert .
$$

In order to prove \eqref{eq-lem-anticoncentration-claim2} for every $k\in\{1,\ldots, p\}$,
we construct  a relation $R_k$ between the sets $\widetilde{\Gamma}_k$
and $\widetilde{\Gamma}_{k-1}$.

Let $\nrangr\in \widetilde{\Gamma}_k$. Note that
\begin{equation}\label{eq-lem-anticoncentration-Sk-1}
\big\vert \sigma([2\vert I\vert]\cup\{i_1,\ldots,i_s\})\cup \innbr\big( J\big)\big\vert \leq 2\vert I\vert +s+d\vert J\vert \leq \frac{9d\vert J\vert}{8}.
\end{equation}
By \eqref{eq-choiceI-lem-anticoncentration}, we get
\begin{equation}\label{eq-lem-anticoncentration-Sk-2}
\big\vert I^c\cap J^c \setminus \outnbr(\sigma(i))\big\vert \geq
n- \frac{d|J|}{32} - \frac{n}{4d} - d   \geq    \frac{27n}{32} - \frac{d|J|}{32} .
\end{equation}

Denote
$$
   S_k:= \edg\big( \sigma([2\vert I\vert]^c\setminus \{i_1,\ldots,i_s\})\setminus \innbr\big(J\big), I^c\cap J^c \setminus \outnbr(\sigma(i_s))\big).
$$
Using (\ref{eq-lem-anticoncentration-Sk-1}), (\ref{eq-lem-anticoncentration-Sk-2}) we observe that
\begin{equation}\label{eq-lem-anticoncentration-size-Sk}
 \vert S_k\vert \geq d \left(\frac{27n}{32} - \frac{d|J|}{32}-\frac{9d\vert J\vert}{8}\right)
\geq \frac{nd}{2}.
\end{equation}

We say that $(\nrangr,\nrangr')\in R_k$ for some $\nrangr'\in \widetilde{\Gamma}_{k-1}$ if $\nrangr'$
can be obtained from $\nrangr$ in the following way.
Let $(\sigma(i_s), j_1)$ be one of the $k$ edges in $\edg(\sigma(i_s),J)$.
Destroy an edge $(v,j)\in S_k$. Since $j\not\in \outnbr\big(\sigma(i_s)\big)$,
then we can create the edge $(\sigma(i_s), j)$.
Since $v\not\in \innbr(j_1)$, then we can destroy the edge $(\sigma(i_s), j_1)$ and
create the edge $(v,j_1)$, thus completing the simple switching.
Therefore by \eqref{eq-lem-anticoncentration-size-Sk} we get
   \begin{equation*}\label{eq-imageRk-anticoncentration}
  \vert  R_k(\nrangr)\vert\geq \frac{knd}{2}.
  \end{equation*}
Note that the above transformation of $\nrangr$ does not decrease
$\vert \innbr(J)\vert$ which guarantees that $\nrangr'\in  \Lambda \big(\frac{1}{2}, J\big)$.

Now we estimate the cardinalities of preimages.  Let $\nrangr'\in R_k(\widetilde{\Gamma}_k)$.
 In order to reconstruct a possible $\nrangr$ for which $(\nrangr,\nrangr')\in  R_k$,
destroy an edge $(v,j_1)$ from $\edgpr(\sigma([2\vert I\vert]^c\setminus\{i_1,\ldots,i_s\}), J)$
to create the edge $(\sigma(i_s), j_1)$ for $j_1\in J$. There are at most
$d\vert J\vert $ such choices. To complete the simple switching,
we destroy an edge $(\sigma(i_s), j_2)$ in $\edgpr(\sigma(i_s), I^c\cap J^c)$ and create the
edge $(v, j_2)$. There are at most $d$ possible choices here.
Therefore
 \begin{equation*}\label{eq-preimageRk-anticoncentration}
  \vert  R_k^{-1}(\nrangr')\vert\leq d^2\vert J\vert .
  \end{equation*}
Claim~\ref{multi-al} implies the inequality \eqref{eq-lem-anticoncentration-claim2}, and completes the proof.
%
\qed

\section{Adjacency matrices of random digraphs}
\label{fromgraph}

In this section we continue to study density properties of random $d$-regular directed (rrd)  graphs. We interpret results obtained in the previous section in terms of adjacency matrices  and provide consequences of the anti-concentration property, Theorem~\ref{th-anticoncentration}, needed to investigate the invertibility of adjacency matrices.

\subsection{Notation}
\label{nota}

For $1\leq d\leq n$ we denote by $\mathcal{M}_{n,d}$ the set of
$n\times n$ matrices with $0/1$-entries and such that every row
and every column  has  exactly $d$ ones. By a random matrix on
$\mathcal{M}_{n,d}$ we understand a matrix uniformly distributed on
$\mathcal{M}_{n,d}$, in other words the probability on $\mathcal{M}_{n,d}$
is given by the normalized counting measure.
Whenever it is clear from the context, we
usually use the same letter $M$ for an element of $\Mc$ and for a random matrix.

For $I\subset [n]$ by $P_I$ we
denote the orthogonal projection on the coordinate subspace $\R^I$ and
$I^c:=[n]\setminus I$. For a  matrix $M\in \Mc$ we say that a non-zero
vector $x$ is a null-vector of $M$ if either $Mx=0$ (a right null-vector) or
$x^T M=0$ (a left null vector).

Let $M=\{\mu_{ij}\}\in \Mc$. The $i$'th row of $M$ is denoted by $R_i=R_i(M)$
and the $i$'th column by $X_i=X_i(M)$, respectively.
For $j\leq n$, we denote $\supp\, X_j=\{i\leq n\,:\, \mu_{ij}=1\}$ and
for every subset  $J\subset [n]$ we let
\begin{equation*}
 S_J:=\bigcup_{j\in J}\,  \supp\, X_j,
\end{equation*}
Clearly,  $|J|\leq |S_J|\leq d|J|$ and $ n-d|J|\leq |(S_J)^c|\leq n-|J|$.

For $x\in \R^n$ we denote its coordinates by $x_i$, $i\leq n$, its  $\ell_{\infty}$-norm by $\|x\|_{\infty}=\max_i |x_i|$
and for a linear operator $U$ from $X$ to $Y$  by $\|U\, :\, X\to Y\|$ we denote
its operator norm.

\subsection{Maximizing columns support}
\label{section-max-support}

In this section we reformulate Theorem~\ref{th-graph-sparse}
in terms of adjacency matrices. It corresponds to bounding from below
the number of rows which are non-zero on a given set of columns. More precisely,
for every subset  $J\subset [n]$ we have $|S_J|\leq d|J|$.
We prove that for almost all matrices in $\Mc$, this inequality is close to being sharp
whenever $J$ is of the appropriate size (less than some proportion of $n/d$).
This means that $S_J$ is of almost maximal size.

\begin{theorem}
\label{graph1}
Let $8\leq d\leq n$ and  $\varepsilon\in (0,1)$ satisfy
$$
\varepsilon^2 \geq \frac{\max\{8, \ln d\}}{d}.
$$
 Define
\begin{equation*}
\Omega_\varepsilon=\Big\{M\in \Mc\,:\,\forall J\subset[n],\,
|J|\leq \frac{c_0 \eps n}{d}, \, \, \text{ one has }\, \,|S_J|\geq
(1-\varepsilon) d|J|\Big\},
\label{Oo}
\end{equation*}
where $c_0$ is a sufficiently small absolute positive constant.
Then
\begin{equation*}
\label{POo}
\mathbb{P}(\Omega_\varepsilon)\geq 1-\exp\left(-\frac{\eps^2 d}{8}
\ln\left(\frac{ec_0\eps n}{d}\right)\right).
\end{equation*}
\end{theorem}

\begin{Remark}
In fact Theorem~\ref{th-graph-sparse} gives slightly more, namely
the corresponding estimates when $|J|=k$ for a fixed
 $k\leq c_0  \varepsilon n/{d}$.
However we don't use it below.
\end{Remark}

The following proposition is a direct consequence of Lemma~\ref{lem1-graph-spread}
(applied with $2\eps$ instead of $\eps$ and with $p=k=0$). It shows
that for a big proportion of matrices in $\Mc$, every two rows have almost disjoint supports.

\begin{prop}\label{lem-disjoint-rows}
Let $\varepsilon\in (0,1)$ and $8\leq d\leq \eps n/6$.
Define
\begin{equation*}
\label{O2e}
\Omega^2_{\varepsilon}=\Big\{M\in \Mc:\,\forall i,j \in [n]\,
 \quad |\supp(R_i+R_j)|\geq 2(1-\varepsilon)d\Big\}.
\end{equation*}
 Then
 \begin{equation*}
 \mathbb{P}(\Omega^2_{\varepsilon})\geq 1- \frac{n^2}{2}\, \left(\frac{e d}{\varepsilon n}\right)^{\varepsilon d}.
\label{PO2e}
\end{equation*}
\end{prop}


\subsection{Large zero minors}
\label{zero-minors}

In this section we reformulate Theorem~\ref{th-graph-spread}
in terms of adjacency matrices.
It states that almost all matrices in $\Mc$ do not contain
large zero minors.

Given $0\leq \alpha,\beta\leq 1$ we define
\begin{align}
\eoo(\alpha,\beta)=\{M\in\Mc\,\, :\, \, &\exists I,J\subset[n]\quad\text{such that}\quad
                     |I|\geq \alpha n,\,|J|\geq \beta n, \notag
\\
&\text{ and }\quad \forall i\in I\, \forall j\in J\quad \mu_{ij}=0\}.
    \label{eo-def}
\end{align}
In other terms, the elements of $\eoo(\alpha, \beta)$ are the matrices in $\Mc$ having a
zero submatrix of size at least $\alpha n\times \beta n$. Theorem~\ref{th-graph-spread}, reformulated below,  shows that this set is small
whenever $\alpha$ and $\beta$ are not very small.

\begin{theorem}
\label{graph2}
There exist absolute positive constants $c,  C$
such that the following holds. Let $2\leq d\leq  n/24$ and
 $0<\alpha\leq \beta\leq 1/4$.  Assume that
$$
\alpha \geq  \frac{C  \ln (e/\beta)}{d}.
$$
Then
$$
\mathbb{P}\left(\eoo(\alpha,\beta)\right) \leq
\exp\left(-c \alpha \beta dn\right).
$$
\end{theorem}

\begin{Remark}
We usually apply this theorem with the following choice of parameters:
$\alpha = p/(2q)$, $\beta = p/2$, where $q= c_1 p^2 d$ for a sufficiently small absolute
positive constant $c_1$.
Then we have
\begin{equation}
\label{remzero}
\mathbb{P}\left(\eoo\left(\frac{p}{2q},\frac{p}{2}\right)\right) \leq  \exp(-c_2 n).
\end{equation}
\end{Remark}

 We will also need the following simple lemma.

\begin{lemma}
\label{graph3}
Let $1\leq d\leq n$ and $0<\alpha, \beta<1$. Let
\begin{align*}
   \Omega_{\alpha,\beta}=\Big\{M\in \Mc\,:\,\,\,
   \forall J,
   \, \,   |J|\geq\beta n, \, \,
 \mbox{ one has } \, \,
  |\{i\,:\,|\supp\,
     R_i\cap J|\geq\beta/2\alpha\}|\geq(1-\alpha)n\Big\}.
\end{align*}
Then provided that $\alpha n$ is an integer, we have
$$
  \left(\eoo(\alpha,\beta/2)\right)^c \subset  \Omega_{\alpha,\beta}.
$$
\end{lemma}

\medskip

\noindent
{\bf Proof. }
Let $M\in \Omega_{\alpha,\beta}^c$. Then there exist $J\subset [n]$ with $|J|\geq \beta n$
and $I\subset [n]$ with $|I|= \alpha n$ such that
$$
  \forall i\in I\, \,  \quad  |\supp\, R_i\cap J|<\beta/2\alpha .
$$
This shows that the minor $\{\mu_{ij}\,:\, i\in I,\, j\in J\}$ has strictly less than
$\beta n/2$ ones, which means that at least $\beta n/2$ columns of this minor
are zero-columns. Thus
$$
  \exists I\subset [n], \,|I|= \alpha n, \,\,\, \, \exists J_0 \subset [n], \,|J_0|\geq \beta n/2,
  \,\,\, \,   \forall i\in I,\, \,   \forall  j\in J_0\,  :\,  \mu_{ij}=0.
$$
In other words, there is a zero minor of size $\alpha n\times \beta n/2$.
This proves the lemma.
\qed

\subsection{An anti-concentration property for adjacency matrices}
\label{acvec}


For every $F\subset [n]\times[n]$ and  $I\subset[n]$, let
$$
  \Mc(I,F)=\left\{\M=\{\mu_{ij}\}\in\Mc \, :\, \mu_{ij}=1 \, \, \mbox{ if and only if }
  \, \, j\in I,\, (i, j) \in F\right\}.
$$
Thus matrices in $\Mc(I,F)$ have the same columns indexed by $I$ and the places of ones in
these columns are given by $F\cap ([n]\times I)$. Of course this set
can be empty.

 For every $\M\in \Mc$, $J\subset [n]$ and $i\leq n$, we define $\delta_i^{J}\in \{0,1\}$ as follows
$$
\delta_i^J=\delta_i^J(\M) :=\left\{
        \begin{array}{ll}
                1  & \mbox{if } \supp R_i\cap J\ne \emptyset ,\\
                0 & \mbox{otherwise.}
        \end{array}
\right.
$$
We also denote $\delta^J:=(\delta_1^J,\ldots, \delta_n^J)\in \{0,1\}^n$.
The quantity $\delta^J$ indicates the rows whose supports intersect with $J$, i.e. the rows that have  at least  one $1$ in columns indexed by $J$.
The following  is a reformulation of Theorem~\ref{th-anticoncentration}, concerning
the anti-concentration property of  graphs, in terms of adjacency matrices.

\begin{theorem}\label{th-anticoncentration-reformulation}
There are absolute positive constants $c, \tilde c$ such that the following holds.
Let $32\leq d\leq cn$ and $I, J$ be disjoint subsets of $[n]$ such that
\begin{equation}\label{eq-choiceI-th-anticoncentration-reformulation}
 \vert I\vert \leq \frac{d\vert J\vert}{32} \quad \text{ \ and \ }\quad
 8\leq \vert J \vert  \leq \frac{8 c n}{d}.
\end{equation}
Let $F\subset [n]\times[n]$ be such that $\Mc(I,F)\neq\emptyset$ and  $v\in \{0,1\}^n$. Then
$$
\P\{ \delta^J =v  \mid  \Mc(I,F)\}\leq 2\exp\left(-\tilde c d\vert J\vert \ln\left(\frac{n}{ d\vert
    J\vert }\right)\right).
$$
\end{theorem}

This theorem has the following consequence.

\begin{prop}\label{th-anticoncentration-matrix part}
There are absolute positive constants $c, c'$ such that the following holds. Let $32\leq d\leq cn$, 
$\lambda\in \R$, $a>0$, and $I, J,J_\lambda$ be a partition of $[n]$ satisfying
 (\ref{eq-choiceI-th-anticoncentration-reformulation}).
 Let $q\leq n/2$ be such that
\begin{equation}\label{eq-choice-q-anticoncentration-matrix}
2^{q+1} \leq \exp\left(  c' d\vert J\vert\ln\left(\frac{n}{d\vert  J\vert}\right)\right)
\end{equation}
and $y$ be a vector in $\R^n$ satisfying
\begin{equation}
\label{eq-condition-y-th-matrix part}
\forall \ell\in J_\lambda\ \ y_\ell =\lambda  \text{\quad and \quad } \forall j\in J\ \ y_j-\lambda \geq 2a.
\end{equation}
Then for every $S\subset [n]$ with $\vert S\vert \geq n-q$, one has
\begin{equation}\label{PS}
\P\{ \Vert P_S\M y\Vert_\infty < a\} \leq \exp\left( -c'd\vert J\vert\ln\left(\frac{n}{d\vert J\vert}\right)\right).
\end{equation}
\end{prop}

\begin{Remark}\label{Remark 1 after th-anticoncentration-matrix part}
The above statement with essentially the same proof holds when
\eqref{eq-condition-y-th-matrix part} is replaced by  \begin{equation*}
\forall \ell\in J_\lambda\ \ y_\ell =\lambda  \text{\quad and \quad } \forall j\in J\ \ \lambda - y_j \geq 2a.\end{equation*}
\end{Remark}


To prove Proposition~\ref{th-anticoncentration-matrix part} we need the following lemma.

\begin{lemma}\label{lem-from matrix to anticoncentration}
Let $\lambda\in \R$, $a>0$, and $I, J,J_\lambda$ be a partition of $[n]$ satisfying (\ref{eq-choiceI-th-anticoncentration-reformulation}).
Let $y$ be a vector in $\R^n$ satisfying (\ref{eq-condition-y-th-matrix part}).
Then for every $i\leq n$ and every  $F\subset [n]\times[n]$ there exists $v_i\in\{0,1\}$ such that
\begin{equation}\label{eq-lem-from matrix to anticoncentration}\nonumber
\{\M\in \Mc(I,F)\, \mid\, \delta_i^J(M)=v_i\}\subseteq
\{\M\in\Mc(I,F)\, \mid \,  \vert (\M y)_i\vert\geq a\}.
\end{equation}
\end{lemma}

\proof
Fix $i\in [n]$ and $F\subset [n]\times[n]$. We argue by contradiction.
Assume that the above inclusion
is violated in both cases,  $v_i=0$ and $v_i=1$. Then there exist two matrices
$\M_1, \M_2\in \Mc(I,F)$ such that
$$
 A_1:=\supp R_i^1\cap J\neq \emptyset, \quad \supp R_i^2\cap J= \emptyset, \quad
  \vert (\M_1 y)_i\vert <a \quad \text{ and } \quad \vert (\M_2 y)_i\vert <a,
$$
where $R_i^j=R_i(M^j)$  denotes the $i$-th row of $M_j$, $j=1,2$. Note that since $M_1, M_2\in \Mc(I,F)$ then
$$
  \supp R_i^1\cap I=\supp R_i^2\cap I := A_2.
$$
Let $s_1:= |A_1|$ and  $s_2:= |A_2|$. Using (\ref{eq-condition-y-th-matrix part}), we observe
$$
(\M_1y)_i = \sum_{j\in A_1} y_j + \sum_{j\in A_2} y_j + \lambda (d-s_1-s_2) \quad \text{ and }\quad (\M_2 y)_i =\sum_{j\in A_2} y_j + \lambda (d-s_2).
$$
Therefore,
$$
 (\M_1y)_i -(\M_2y)_i=  \sum_{j\in A_1} (y_j-\lambda) \geq 2s_1 a\geq 2a,
$$
which is impossible as  $ \vert (\M_1 y)_i\vert <a \text{ and } \vert (\M_2 y)_i\vert <a$.
\qed

\medskip

\noindent{\bf {Proof of Proposition~\ref{th-anticoncentration-matrix part}.}}$\quad$
Since (\ref{eq-condition-y-th-matrix part}) and (\ref{PS}) are homogeneous in $y$, without loss of generality we assume that $a=1$.

Fix $S\subset [n]$ with $|S|\geq n-q$. Let $\cal{F}$ be the set of all
$F\subset [n]\times[n]$ such that $\Mc(I,F)$ is not empty.
Note that
$\{\Mc(I,F)\}_{F\in \mathcal{F}}$ form a partition of $\Mc$. Therefore
it is enough to prove that for every $F\in \mathcal{F}$,
$$
p_0 := \P\{ \Vert P_S\M y\Vert_\infty < 1\mid \Mc(I,F)\} \leq \exp\left( -c'd\vert J\vert\ln\left(\frac{n}{d\vert J\vert}\right)\right).
$$
Fix $F\in \mathcal{F}$. Let
$v_1,\ldots, v_n\in \{0,1\}$ be  given by
Lemma~\ref{lem-from matrix to anticoncentration}.
Note that
\begin{equation*}\label{eq1-proof-th-anticoncentration-matrix part}
 \Vert P_S\M y\Vert_\infty <1\quad\quad \text{ iff }\quad \quad
 \forall i\in S \quad \vert (\M y)_i\vert < 1,
\end{equation*}
therefore if $\|P_S\M y\|_\infty <1$ then
$\{i\, :\,\,  \delta_i^J(M) =v_i\}\subset S^c$.
Thus
$$
  p_0 \leq \P\{ \{i:\ \delta_i^J(M)= v_i\}\subseteq S^c \mid \Mc(I,F)\}.
$$
Now for every $K\subset [n]$, define $v^K\in \{0,1\}^n$ by
$$
v_i^K:=\left\{
        \begin{array}{ll}
                v_i  & \mbox{if } i\in K, \\
                1-v_i & \mbox{otherwise.}
        \end{array}
\right.
$$
Since $m:=\vert S^c\vert \leq q$, by
Theorem~\ref{th-anticoncentration-reformulation} we obtain
\begin{align*}
p_0
&\leq \sum_{\ell =0}^{m}  \P\{ \exists K\subset S^c:\ \vert K\vert =\ell \text{ and }\delta^J(M)= v^K \mid \Mc(I,F)\}\nonumber \\
& \leq \sum_{\ell =0}^{m} {m \choose \ell} \max_{\vert K\vert =\ell} \P\{\delta^J(M)= v^K \mid \Mc(I,F)\}
\leq  2^{q+1}\exp\left( -\tilde c d\vert J\vert\ln\left(\frac{n}{d\vert J\vert}\right)\right).
\end{align*}
Taking $c'=\tilde c/2$ and using (\ref{eq-choice-q-anticoncentration-matrix}) we complete the proof.
\qed

\section{Invertibility of adjacency matrices}
\label{s:matrices}

In this section we investigate the invertibility of adjacency
matrices $M\in\Mc$ of random $d$-regular directed graphs and prove Theorem~A.

\subsection{Almost constant null-vectors}
\label{s:AC}

 We say that  a non-zero vector  is ``almost constant" if for some
$0<p<1/2$
at least $(1-p)n$ of its coordinates
are equal to each other.
Formally,
for $0 < p <1/2$ consider the following set of vectors
\begin{equation}
\label{AC}
   AC(p)=\{x\in\R^n\setminus\{0\}\,:\,\exists
   \lam_x\in\R\quad |\{i\ :\;x_i=\lam_x\}|\geq
   (1-p)n\}.
\end{equation}
In this section we estimate the probability of the event
\begin{equation}
\label{eAC}
    \eac:=\{M\in\Mc\,:\, \forall x\in AC(p)\quad\M x\neq 0
 \quad\text{and}\quad x^T\M\neq0\},
\end{equation}
which relates almost constant vectors to null vectors of $M$.
We show that that this probability is close to one, in other
words we show that
with high probability
a matrix $M\in\Mc$ cannot have almost constant null vectors.
This will be used in the proof of the main theorem allowing one
to restrict the proof to the event $\eac$. More precisely,
we prove the following theorem.

\begin{theorem}
\label{t:AC}
There are absolute positive constants $C, c$ such that for $C\le d\le cn$ and $p\le c/\ln d$ one has
\begin{equation}
     \p\big(\eac\big)\geq 1-\left(\frac{Cd}{n}\right)^{cd}.
  \label{pAC}
\end{equation}
\end{theorem}

\medskip


We start with some comments on the structure of
almost constant vectors. Since $p<1/2$, if
$x\in AC(p)$ then only one real number $\lam_x$ satisfies (\ref{AC}).
For every $x\in AC(p)$ we fix such $\lam_x \in \R$. We set
 \begin{equation*}
   AC^+(p)=\{x\in AC(p)\,:\,\lam_x\geq 0\}.
\end{equation*}
Note that $\lam_{-x} = - \lam_{x} $, therefore
$$
 \eac=\{M\in\Mc\,:\, \forall x\in AC^+(p)\quad\M x\neq 0
 \quad\text{and}\quad x^T\M\neq0\}.
$$
Moreover,
since $(x^TM)^T=M^Tx$ and $M^T$ has the same distribution as $M$ then
\begin{align*}
 \P(\{M\in\Mc\,:\, \forall x\in AC^+(p)\quad\M x\neq 0
    \})=\P(\{M\in\Mc\,:\, \forall x\in AC^+(p)\quad x^T\M \neq 0
    \}).
\end{align*}
Therefore it is enough to consider the event
\begin{align*}
  \e^{AC+}(p)=\{M\in\Mc\,:\, \forall x\in AC^+(p)\quad\M x\neq 0
    \}\notag\\
 \end{align*}
 and to prove that
 \begin{equation*}
     \p\big(\e^{AC+}(p)\big)\geq 1-\frac{1}{2}\left(\frac{Cd}{n}\right)^{cd}.
  \label{pAC+}
\end{equation*}
To this end we  split $AC^{+}(p)$ into two complementary sets and
treat them separately in two lemmas.

\bigskip

For a vector $x = (x_i) \in\R^{n}$ we define the rearrangement
$x^\star=(x_i^\star)_i $ as follows: $x_i^\star = x_{\pi(i)}$, where
$\pi:[n]\to[n]$ is a permutation of
$[n]$ such that  $(|x_i^\star|)_i$ is a decreasing sequence,
that is, $\vert x_{\pi(1)}\vert \ge \vert x_{\pi(2)}\vert \ge \ldots \geq \vert x_{\pi(n)}\vert$.
Contrary to the usual decreasing rearrangement of absolute values
of a sequence, here values $x_i^\star$ can be negative.

In the proof, we choose appropriately a positive integer $m_0$ and consider
a certain subset of $AC^{+}(p)$. For a vector $x$ in this subset
we ``ignore" its first $m_0$ coordinates $x^\star _i$, i.e. we consider $P_I x^\star$ with $I=[m_0]^c$.
Then we show that this vector can be split into a sum of two vectors with disjoint supports and such that the second vector has equal coordinates
on its support. To approximate such vectors in $\ell _\infty$-metric  we construct a net in the following way.

Let $\eta >0$ be a reciprocal of an integer.
 For every $H\subset [n]$ of cardinality $k_1:=pn - m_0$ (we choose $p$ so that $pn$ is an integer) fix an $\eta$-net
$\Delta_H$   in the cube $P_H ([-1, 1]^n)$.
Such $\Delta_H$ can be chosen with $|\Delta_H|\le(1/\eta)^{ k_1}$.
Given $L\subset [n]$ of cardinality $k_2:=(1-p)n$,  consider  the one-dimensional space generated by the vector $v_L$ with $\supp v_L=L$
and all coordinates on $L$ equal to one. Fix an $\eta$-net $\Lambda_L$ in the segment $[-v_L,\,v_L]$. Clearly, $\Lambda_L$ can be chosen with $|\Lambda_L|= 1/\eta $. Note also that for every $z\in \Lambda_L$ one has $\supp z=L$ and $z_i=z_j$ whenever $i, j\in L$, that is $z\in AC(p)$ and $z_i=\lam_z$ for $i\in \Lambda$. Given disjoint subsets $H$, $L$ of $[n]$ of cardinalities $k_1$ and $k_2$
respectively, consider $\Delta_H\oplus\Lambda_L=\{w+z \, :\, w\in \Delta_H, \, z\in \Lambda_L \}$. Then
$$
   \Delta_H\oplus\Lambda_L  \subset AC(p) \quad \quad \mbox{ and }\quad \quad  \left| \Delta_H\oplus\Lambda_L \right|\leq (1/\eta)^{ k_1+1}\leq(1/\eta)^{pn}.
$$
Finally we observe that the vector $P_Ix^\star$ can be approximated by the vectors in the union of $\Delta_H\oplus\Lambda_L$ over all such choices of $H$ and $L$.

In fact we will use only a subset of this union.
Fix a parameter $a>0$ and a positive integer $r$. For
$H$, $L$ as above 
consider
$$
 \Gamma(H, L)=\Gamma_{a,r}(H, L):=\{y\in  \Delta_H \oplus \Lambda_L\,:\, \exists
 J\subset H, \, |J|=r,\,\,
 \text{ such that }\, \,  \forall i\in J \,\, \, y_i-\lambda_y\geq 2a\}.
$$
Clearly, $\vert \Gamma(H,L)\vert \leq (1/\eta)^{pn}$.
Finally, set
\begin{align*}
\mathcal{N} = \mathcal{N}_{a,r}:= \bigcup_{|L|=k_2,|H|=k_1}
\Gamma(H, L),
\end{align*}
where the union is taken over all disjoints subsets $H$ and $L$ of $[n]$ of cardinalities $k_1$ and $k_2$ correspondingly.
Then
\begin{equation}\label{eq-size-subnet}
\vert {\mathcal{N} }\vert
\leq {n\choose k_2}\, {n-k_2 \choose k_1}\, \left(\frac{1}{\eta}\right)^{pn}
\leq {n \choose pn}\, {pn \choose m_0}\, \left(\frac{1}{\eta}\right)^{pn}
\leq \left(\frac{2 e}{\eta  p}\right)^{pn}.
\end{equation}

\smallskip

We are ready now to prove two lemmas needed for Theorem~\ref{t:AC}.
In both of them we use the following set associated with  $x\in AC^+(p)$
and a given $m_0$,
$$
     J_x =J_x(m_0) :=\{i> m_0\;:\;|x_i^\star-\lambda_x| \geq 1/(2d)\}.
$$

\begin{lemma}
\label{l:T2}
There are absolute positive constants $c$ and $c_1$ such that the
following holds.
Let $32\leq d\leq cn$, $m_0\geq 1$, and $r\geq 8$ be integers such that
 $1\leq 2 c_1 r \ln(n/(dr))$. Let $p\in (0, 1/2)$ be such that $pn$ is an integer.
 Assume that
\begin{equation}
\label{eq-ac-condition-p-lem1}
m_0\leq 2 c_1 r \ln\left(\frac{n}{dr}\right),\quad
r\leq \frac{8 c n}{d}, \quad
p\leq \frac{d r}{32 n}, \ \mbox{ and } \
p\left(\ln(e/p)+ \ln(18 d^2)\right)\leq \frac{c_1 d r}{n}\ln \left(\frac{n}{dr}\right).
\end{equation}
Consider the following subset of almost constant vectors
\begin{align*}
T_1=\{x\in AC^+(p)\,:\,\,&|x^\star_{m_0}|=1\quad\text{and}\quad
                          |J_x|\geq 2r\}\label{2b}
\end{align*}
and the corresponding event
\begin{align*}
\e_{T_1}=\{M\in\Mc\,:\, \forall x\in T_1\quad\M x\neq 0\}.
\end{align*}
Then
\begin{equation*}
\mathbb{P}(\e_{T_1})\geq 1-2\exp\left(- c_1dr \ln\left(\frac{n}{d r}\right)\right).
\end{equation*}
\end{lemma}

\begin{Remark}
We apply this lemma with $r=c_2 n/d$,
$m_0 = c_3 n/d$, so that the probability is exponentially (in $n$)
close to one.
\end{Remark}

\begin{Remark}\label{Remark 2 after lT2}
In fact we show a stronger estimate which
could be of independent interest, namely
\begin{align*}
\p \big(\{M\in\Mc\,:\, \exists x\in T_1
\text{ \ such that \ }
& \| Mx \|_{\infty} < 1/(8d)  \}\big)
\le
2\exp\left(- c_1dr \ln\left(\frac{n}{d r}\right)\right). 
\end{align*}
\end{Remark}

 \newcommand{\Quk}{P_{S_{K}^c}}

\noindent
{\bf Proof of Lemma~\ref{l:T2}. \ }
We prove a stronger bound from Remark~\ref{Remark 2 after lT2}.
We start by few general comments on the
strategy behind the proof.
 By the construction of $T_1$, for  $x=(x_i)_i\in T_1$ we have
\begin{align*}
\max\Big\{|\{i> m_0\,:\,x_i^\star-\lambda_x\geq 1/(2d)\}|,\;|\{i>
  m_0\,:\,\lambda_x-x_i^\star\geq 1/(2d)\}|\Big\}\geq r.
\end{align*}
Therefore denoting
$$
   T_1^+:=\{x\in T_1\, :\, |\{i> m_0\,:\,x_i^\star-\lambda_x\geq 1/(2d)\}| \geq r\}
$$
and
$$
T_1^-:=\{x\in T_1\, :\, |\{i> m_0\,:\,\lambda_x-x_i^\star \geq 1/(2d)\}| \geq r\},
$$
we have  $T_1\subseteq T_1^+\cup T_1^-$.
Thus it is sufficient to show that
\begin{align}
p_0:= \p \big(\{M\in\Mc\,:\, \exists x\in T_1^+ \,\,\,\,
& \| Mx \|_{\infty} < 1/(8d) \}\big)
\le
\exp\left(- c_1dr \ln\left(\frac{n}{d r}\right)\right)
\label{more+}
\end{align}
and similarly for $T_1^-$. Below we prove (\ref{more+}) only.
Its counterpart for $T_1^-$ follows the same lines, one just needs to apply Proposition~\ref{th-anticoncentration-matrix part}
with Remark~\ref{Remark 1 after th-anticoncentration-matrix part} below
(with a slight modification of the net constructed above).

\smallskip


To prove (\ref{more+}) we first  approximate vectors in $T_1^+$ by elements of the net ${\cal N}$ constructed above. By the union bound, this will reduce (\ref{more+}) to estimates on the net.
Then, applying Proposition~\ref{th-anticoncentration-matrix part},
we obtain a probability bound for a fixed vector from the net.
As usual, the balance between the probability bound and the size of the net plays the crucial role.

\smallskip

Fix two parameters $\eta:=1/(9 d^2)$ and $a= 1/(4d)-\eta$, and take $k_1=pn-m_0$, $k_2=(1-p)n$ as in the construction
of the net $\cal N$ above.
We start by showing how an element of $T_1^+$ is approximated by
an element from $\cal N$.
Let $x\in T_1^+$ and assume for simplicity that
$\vert x_1\vert\ge \vert x_2\vert\ge \ldots \ge \vert x_n\vert$ (that is, $x=x^\star$).
Recall that $\lambda_x$ is the unique real number satisfying
 (\ref{AC}). By the definition of $T_1^+$ it is easy to see
that there exists a partition $J, J_0, I$ of $[n]$ such that
\begin{align}
&|J|= r,
\quad |J_0| = k_2,
\quad |I| =n-r-k_2,
\notag \\
&\forall i\in J_0\quad x_i=\lam_x \notag \
\quad{\text{with} \ }\quad  \lam_x\ge 0, \\
 &\forall j\in J\quad \,
j>m_0 \quad{\text{and} \ }\quad x_j\geq \lam_x+1/(2d).\notag
\end{align}

Since $|x_{m_0}|=1$ and there is $i> m_0$ such that $x_i\geq \lam_x+1/(2d)$, we observe that $\lam_x< 1$. Since for $i\leq m_0$ we have either $x_i\ge 1$ or $x_i\le -1$, then $J_0\cap I_0=\emptyset$, where $I_0=[m_0]$. Note also that $J\cap I_0=\emptyset$, hence $I_0\subset I$. Set $H=J\cup (I\setminus I_0)$ and $L= J_0$. Then $\vert H\vert =k_1$, $\vert L\vert =k_2$, and $A:=I_0^c= H\cup L$.
By the definition of  $\Delta_H$ and $\Lambda_L$ there exist
$y'\in\Delta_H$ and $y''\in\Lambda_L$ such that
$$
    \| P_{H}x -  P_{H}y'\|_{\infty} \le \eta \quad \quad \mbox{ and } \quad \quad\| P_{L} x -  P_{L} y''\|_{\infty} \le \eta.
$$
Therefore $y:=y'+y''\in \Delta_{H}\oplus \Lambda_{L}$ satisfies
$\| P_{A}x -  P_{A}y\|_{\infty} \le \eta$.
Moreover, by the construction of the net $y\in AC(p)$,
$$
   \forall i\in L\quad \ y_i=y''_i=\lambda_y\quad    \text{\quad \ and \quad } \quad \forall i\in J\,\,\, y_i-\lambda_y \geq x_i -\lambda_x -2\eta \geq
  (2d)^{-1}-2\eta=2a.
$$
Thus we showed that for every $x\in T_1^+$ there exist  $H, L \subset [n]$ with $|H|=pn-m_0$, $|L|=(1-p)n$, and $y\in \Gamma (H, L) =     \Gamma_{a,r} (H, L)$ such that $\|P_A x - P_A y\|\leq \eta$. Note also, that given
$H$ and $L$ one can ``reconstruct" $I_0$ as $I_0=[n]\setminus (H\cup L)$.

%
%
Moreover, denoting
\begin{equation*}
 S:=S_{I_0}^c=[n]\setminus \supp\sum_{i\in I_0}X_i.
\end{equation*}
and observing that $P_{S}MP_{I_0}=0$ (indeed, for every $i\in S$ and $j\in I_0$ one has   $\mu_{ij}=0$), we get
\begin{align*}
\Vert P_{S} M y\Vert_{\infty}
&= \Vert P_{S} M x +  P_{S} M (y-x)\Vert_{\infty}
=  \Vert P_{S} M x +  P_{S} M P_{A}(y-x)\Vert_{\infty}\nonumber\\
&\leq  \Vert P_{S} M x\Vert_{\infty}+ \Vert P_{S} M P_{A}(y-x)\Vert_{\infty}
<\|Mx\|_{\infty} + \Vert  M\, : \, \ell_{\infty} \rightarrow \ell_{\infty} \Vert \, \eta
\nonumber\\&
\leq 1/(8d) +\eta d < a,
\end{align*}
provided that $\|Mx\|_{\infty}\leq 1/(8d)$.
Thus, by the union bound, we obtain
$$
  p_0 \leq \sum _{y\in {\mathcal{N}}} \p\big(\{M\in\Mc\,\, :\, \, \Vert P_{S} M y \Vert_{\infty} <a\}\big),
$$
where $S=S(y)=S^c_{I_0}$, $I_0=I_0 (y)=[n]\setminus (H\cup L)$ whenever
$y\in \Gamma (H, L)$.

%
%


Finally we estimate the probabilities in the sum. Let $H, L \subset [n]$ be such that $|H|=pn-m_0$, $|L|=(1-p)n$, and $y\in \Gamma (H, L)$,  $J$ be from the definition $\Gamma (H, L)$ and $S$ be as above. Let $I=
[n]\setminus (J\cup L)$. Then $I$, $J$, $L$ form a partition of $[n]$ with $|J|=r$ and $|I|= pn-r$.  By assumptions of the lemma, this partition
satisfies (\ref{eq-choiceI-th-anticoncentration-reformulation}).
Note also that assumptions on $m_0$ and $r$ imply $m_0 d< n$, hence  $\vert S\vert \geq n-m_0 d>0$, and
 (\ref{eq-choice-q-anticoncentration-matrix}) is satisfied with  $q:=m_0 d$ (with $c_1=c'/2$).
By the definition of $\Gamma (H, L)$ the vector $y$ satisfies
$$
  \forall i\in L \quad y_i=\lambda_y \quad \ \text{ \ and \ } \quad\forall i\in J\, \, \, y_i-\lambda_y\geq 2a.
$$
Applying Proposition~\ref{th-anticoncentration-matrix part}
with the partition $\{I, J, L\}$, the vector $y$, and the
set $S$, we obtain
\begin{equation*}
\p\big(\{M\in\Mc\,:\, \Vert P_{S} M y\Vert_{\infty} <a\}\big)\leq
\exp\left(- 2c_1dr \ln\left(\frac{n}{d r}\right)\right).
\end{equation*}
Since $\eta=1/(9d^2)$, by \eqref{eq-size-subnet} and \eqref{eq-ac-condition-p-lem1}, this implies
\begin{align*}
p_0\leq |{\mathcal{N}}| \, \exp\left(- 2c_1dr \ln\left(\frac{n}{d r}\right)\right) \leq
 \left(\frac{18 d^2 e}{p}\right)^{pn}\, \exp\left(- 2c_1dr \ln\left(\frac{n}{d r}\right)\right)
\leq \exp\left(- c_1dr \ln\left(\frac{n}{d r}\right)\right),
\end{align*}
which completes the proof.
\qed


\smallskip

In the next lemma we prove an analogous statement for the set which is complementary to $T_1$. Recall that $\Omega_{\varepsilon}$ was introduced
in Theorem~\ref{graph1}   and let $c_0$ be the same constant  as  in that theorem.

\begin{lemma}
\label{l:T3}
Let $\varepsilon\in (0,1/4)$. Let $m_0, m_1, r$ be positive integers such that
$$
  m_1=m_0+2r < \min\{ m_0/(2\eps), \, c_0 \eps n/d\}.
$$
Consider the following subset of almost constant vectors
\begin{align*}
T_2=\{x\in AC^+(p)\,:\,\, {\text{either} \ }\quad|x^\star_{m_0}|=0\, \text{ \ or \ } \, (|x^\star_{m_0}|=1\quad\text{and}\quad |J_x|< 2r\, )\}.
\end{align*}
Then
\begin{align*}
\Omega_{\varepsilon}\subset \e_{T_2}:=\{M\in\Mc\,:\, \forall x\in T_2\quad\M x\neq 0\}.
\end{align*}
\end{lemma}

To prove the lemma we need the following simple observation.

\begin{claim}\label{sobs} Let $\eps\in (0, 1/2)$, $1\leq d\leq n$, and
 $1\leq m\leq c_0\eps n/d$, where $c_0$ is the constant from the definition
of $\Omega_{\varepsilon}$.
Let $M\in \Omega_{\varepsilon}$ and $I$ be the set of indices corresponding to rows having exactly one $1$ in columns indexed by $[m]$, i.e.
 \begin{equation*}
 I=\{i\in S_{[m]}\, :\, |\supp R_i\cap [m]|=1\}.
\end{equation*}
Then $\vert I\vert\geq (1-2\eps) d m >0.$
\end{claim}

\noindent
{\bf Proof.  }
Since $M\in \Omega_{\varepsilon}$,
\begin{equation*}
|S_{[m]}|\geq (1-\varepsilon)dm.
\label{Smx}
\end{equation*}
 Since rows $R_i$, $i\in I$, have exactly one $1$ on $[m]$, while
other rows indexed by $S_{[m]}$ have at least two ones on $[m]$,
we observe
$$
\vert I\vert +2(|S_{[m]}|-\vert I\vert )\leq dm .
$$
This implies the desired result.
\qed

\noindent
{\bf Proof of Lemma~\ref{l:T3}.  }
Let $M\in\Omega_{\varepsilon}$ and $x\in T_2$.
For simplicity  assume that $x=x^\star$, so that
$|x_1|\geq|x_2|\geq...\geq|x_n|$.
Our proof consists of the following three cases.

\vspace{2ex}
\noindent
{\bf Case 1: $\vert x_{m_0}\vert =0$. \, \, }
Let $m_x\geq 1$ be the largest integer such that $\vert x_{m_x}\vert \neq 0$. Clearly $m_x< m_0$. Let $I_x$ be the set of indices corresponding to rows having exactly one $1$ in columns indexed by $[m_x]$. By Claim~\ref{sobs}, $I_x\ne \emptyset$.
Thus there  exists a row $R_i$, $i\in I_x$, and a unique $j\in [m_x]$ such that
$\mu_{ij}=1$. This implies
$$
  (Mx)_{i} =  \langle R_i,x\rangle =  x_j\ne 0.
$$

\vspace{2ex}
\noindent
{\bf Case 2: $\vert x_{m_0}\vert=1$, $|J_x|<2 r$, and $\lambda_x < 1/({2d})$. \, \, } In this case the cardinality of the set
$$
  \{i\geq m_0\, :\, \vert x_i\vert \geq \lambda_x +1/(2d)\}
$$
is  less than $2r$. Since $m_1-m_0 = 2r$, we have
$ \vert x_{m_1}\vert <\lambda_x +  1/(2d) < 1/d.$

Let $J_j=[m_j]$, $j=0,1$. We first show that there is a row $R_{i}$ such that
\begin{equation}
|\supp R_{i}\cap J_{0}|=1\quad \mbox{ and } \quad |\supp R_{i}\cap(J_{1}\setminus J_{0})|=0.
\label{Rik}
\end{equation}

Let $I$ be the set of indices corresponding to rows having exactly one $1$ in $J_1$. By Claim~\ref{sobs}, $|I|\ge (1-2\eps) d m_1$.
Since the number of nonzero rows on $J_1\setminus J_{0}$ is at most $d\vert J_1\setminus J_{0}\vert$,  the number of rows satisfying (\ref{Rik}) is at least
$$
  (1-2\eps) d m_1- d(m_1-m_{0})= d(m_{0}- 2\eps m_{1}) >0,
$$
that is  there exists a row $R_{i}$ satisfying (\ref{Rik}).  Denote
$j_0\in J_{0}$ the only coordinate of $P_{J_1 }R_{i}$ which is equal to one,
i.e. $\mu_{i j_0}=1$ and for every
$j\in J_1\setminus\{j_0\}$, $\, \mu_{i j}=0$. Therefore if $j\in \supp R_{i}$ and $j\neq j_0$
then $j>m_1$ and $\vert x_j\vert \leq \vert x_{m_1}\vert$.
Since  $\vert x_{m_1}\vert < 1/d$, we observe
$$
 |(Mx)_{i}| =  \vert \langle R_{i},x\rangle\vert \geq \vert x_{j_0}\vert -\sum_{\underset{j\neq j_1}{j\in \supp R_{i_1}}} \vert x_j\vert
\geq |x_{m_0}|-(d-1)|x_{m_1}|\geq 1-\frac{d-1}{d}>0.
$$

\vspace{2ex}
\noindent
{\bf Case 3: $\vert x_{m_0}\vert =1$, $|J_x|< 2r$, and $\lambda_x \geq 1/(2d)$. \, \, }
Consider the set
$$
  J:=  \{ i\leq n\,: 0<x_i<\lambda_x+ 1/(2d)\}.
$$
Then $A:=J^c\subset [m_0]\cup J_x$. Thus
$\vert S_{A}\vert \leq (m_0+2r)d <n $.
Therefore, there exists a row $R_j$, $j\notin S_{A}$, such that $\supp R_j\subset J$. This implies
$$
   (M x)_j = \langle R_j, x\rangle =\sum_{j\in J} x_j >0.
$$

Thus in all cases $Mx\ne 0$, which completes the proof.
\qed

\medskip

\noindent{\bf Proof of Theorem \ref{t:AC}.}
Recall that as we mentioned after the theorem it is enough to bound the probability of the event $\e^{AC+}$.

Let $c$, $c_0$, $c_1$ be constants from Lemmas~\ref{l:T2} and \ref{l:T3}.
We choose $\eps>0$ to be small enough constant ($\eps =\min \{1/8, c_1 c_0/4, 32c/c_0\}$ would work), $m_0=\lfloor 2 c_0 \eps ^2 n/d\rfloor$ and
$r=\lfloor m_0/(8\eps )\rfloor\approx c_0 \eps n/(4d)$, so that assumptions
of Lemmas~\ref{l:T2} and \ref{l:T3} on $m_0$ and $r$ are satisfied. Finally,
for a sufficiently small absolute positive constant $c_2$ we choose the biggest $p\leq c_2 /\ln d$ such that $pn$ is an integer. Then assumptions of Lemma~\ref{l:T2} on $p$ are also satisfied (note that it is enough to prove the theorem with the biggest possible $p$). Therefore, applying these lemmas together with Theorem~\ref{graph1},  we have
$$
\mathbb{P}(\e_{T_1})\geq 1- 2\exp(-c_3 n) \text{ \ and \ }
\mathbb{P}(\e_{T_2})\geq \mathbb{P}(\Omega_{\eps})\geq
1-\exp\left(-c_4 d \ln\left(\frac{c_5 n}{d}\right)\right),
$$
where $T_1$, $T_2$ are events from the lemmas and $c_i$'s are absolute positive constants. Since $\e^{AC+}\supseteq \e_{T_1}\cap \e_{T_2}$ we obtain the desired result
by adjusting absolute constants.
\qed

\subsection{Auxiliary results}
\label{aux}

\subsubsection{Simple facts}
\label{facts}

We will need the two following simple facts.

\begin{claim}
\label{c:J}
Let $p\in (0, 1/3]$ and  $x\in\R^n$. Assume that
$$
 \forall \lam\in\R\quad |\{i\,:\,x_i=\lam\}|\leq(1-p)n.
$$
Then there exists $J\subset[n]$ such that
$$
  pn\leq|J|\leq  (1-p)n \quad\text{ and }\quad \forall i\in J\,\, \forall j\notin J \quad
  x_i\neq x_j.
$$
\end{claim}

\begin{Remark}\label{aaaaaa}
We apply this claim twice, once in the following form. Let $m\le n$ and $\ell \leq m/3$. Let $S\subset[n]$, $|S|=m$, and let 
$v\in \R^n$ satisfy $\, \forall \lam\in\R\quad |\{i\in S\,:\,x_i\ne\lam\}|\geq \ell$.
Then there exists $J\subset S$ such that
$$
  \ell \leq|J|\leq  m-\ell \quad\text{ and }\quad \forall i\in J\,\, \forall j\in S\setminus J \quad
  v_i\neq v_j.
$$
\end{Remark}

\noindent
{\bf Proof of Claim~\ref{c:J}.  }
Let $\{\lam_1,...,\lam_k\}$ be the set of all distinct values of
coordinates of $x$.
For  $j\leq k$,
let $I_j=\{i\,:\, x_i=\lam_j\}$ and $m_j=|I_j|$.
Clearly, $m_j\leq(1-p)n$ for every $j$.
By relabeling
assume that  $m_1\geq m_2\geq...\geq m_k$.
If $m_1\geq pn$, choose $J=I_1$. Otherwise set
$J=I_1\cup...\cup I_{\ell}$,
where $\ell$ is the smallest number satisfying
$m:=|J|=m_1+...+m_{\ell}\geq pn$.
Since $m_j\leq m_1<pn$ for $j\leq k$, then
$m< 2pn$, and this implies
$$
 pn\leq|J|<2pn \leq (1-p)n.
$$
\qed

\medskip

Let $A$, $A_1$, ..., $A_m$ be sets such that every $x\in A$ belong
to at least $k$ of sets $A_i$'s. Then we say that $\{A_i\}_i$ forms
a $k$-fold covering of $A$.

The proof of the following fact uses a standard argument in measure theory,
so we omit it.

\begin{claim}
\label{k-fold}
Let $(X, \mu)$ be a measure space. Let $A$, $A_1$, ..., $A_m$ be subsets
of $X$ such that $\{A_i\}_i$ forms a $k$-fold covering of $A$. Then
$$
     k \mu (A) \leq \sum _{i=1}^{m} \mu (A_i).
$$
\end{claim}

\subsubsection{Combinatorial results}
\label{comb}

In this section we prove a Littlewood-Offord type result,
which will be one of key steps in the shuffling procedure.

 Consider a probability space
$$
    \Omega_0 =   \{B\subset [2d]\,:\,|B|=d\}
$$
with the probability given by the normalized counting measure.
For a vector $v\in \R^{2d}$ and $B\in  \Omega_0$ denote
$$
   v_B = \sum _{i\in B} \, v_i .
$$

\begin{prop}
\label{shary}
Let  $1\leq k \leq d$. Let $v\in\R^{2d}$ and $a\in \R$. Assume
there exists $J\subset[2d]$ such that $|J|=k$ and for every
$i\in J$ and every $j\not\in J$ one has $v_i\neq v_j$. Then
$$
    \p (v_B = a) \leq \frac{10}{\sqrt{ k}}.
$$
\end{prop}

\medskip

To prove Proposition~\ref{shary} we need two combinatorial lemmas.
We start with so-called anti-concentration Littlewood-Offord type lemma
(\cite{Er}, see also \cite{Kl}). Usually it is stated for $\pm 1$ Bernoulli random variables, but by shifting and rescaling it holds for any
two-valued random variables, where by a two-valued random variable we mean
a variable that takes two different values, each with probability half.

\begin{lemma}
\label{anti}
Let
$\xi_1$, $\xi_2$, ..., $\xi _m$ be independent two-valued random variables.
Let $x\in \R^m$.
Then
$$
 \sup_{a\in \R}\p\Big(\sum_{i=1}^m\xi_i  x_i=a\Big)\leq |\supp x|^{-1/2}.
$$
\end{lemma}

\medskip

Let $\Pi_{2d}$ be the permutation group with a probability given by the
normalized counting measure and denoted by $\p_{\Pi_{2d}}$. By $\pi$ we
denote a random permutation. The proof of the next lemma is rather
straightforward, we postpone it to the end of the section.

\begin{lemma}
\label{pi}
Let  $1\leq k\leq d$.
Let $x\in\R^{2d}$ and $J \subset[2d]$ be such that
$|J|=k$ and for every $i\in J$ and every $j\not\in J$ one has
$x_i\neq x_j$. For $\pi\in\Pi_{2d}$ let
$$
   E=E(\pi)=\{(x_{\pi(i)},x_{\pi(i+d)})\,: \,  i\leq d,\, \, x_{\pi(i)}\neq x_{\pi(i+d)}\}.
$$
Then
$$
  \p_{\Pi_{2d}}\left( |E|\leq \frac{k}{50} \right)\leq  \left(\frac{k}{1.1 \, d}\right)^{k/3}.
$$
\end{lemma}

\medskip

\noindent{\bf Proof of Proposition~\ref{shary}. } Let $B$ be a (set-valued) random
variable on $\Omega _0$. Let $\delta=(\delta_1,...,\delta_d)$ be a vector of independent
Bernoulli random $0/1$ variables ($\p(\delta_i=1)=1/2$ for $i\leq d$), and $\Omega$
denote the corresponding probability space. Consider a random
(on $\Pi_{2d}\times\Omega$) set of indexes
$$
   A(\delta,\pi)=\{\pi(i)\,:\, \delta_i=1\}\cup\{\pi(i+d)\,:\, \delta_i=0\}\subset [2d].
$$
Note that $|A(\delta,\pi)|=d$.
It is not difficult to see that
for every fixed $\delta$,
$A(\delta,\pi)$ has the same distribution as $B$.
Therefore, $A(\delta,\pi)$ on $\Pi_{2d}\times\Omega$
has the same distribution  as $B$ on $\Omega_0$. Thus,
given $v\in\R^{2d}$, the random variable $v_B$
 has the same distribution as $v_{A(\delta,\pi)}$.
 Now we introduce the following random variables on $\Pi_{2d}\times\Omega$:
$$
 \xi_i=\xi_i^v=\delta_iv_{\pi(i)}+(1-\delta_i)v_{\pi(i+d)}.
$$
Note that $\p(\xi_i=v_{\pi(i)})=\p(\xi_i=v_{\pi(i+d)}) =1/2$ and that
$$
   v_{A(\delta,\pi)}=\sum_{i\in {A(\delta,\pi)}}v_{i}=\sum_{i=1}^d\xi_i.
$$
 Moreover, if we condition on $\pi$, the random variables $\bar{\xi}_i=\xi_i|_\pi$ are independent,
hence we can apply Lemma~\ref{anti}. Denote by $m(\pi)$ the number of  two-valued $\bar{\xi}_i$'s.
Then
$$
  \p_{\Omega}\Big(\sum_{i=1}^d\bar \xi_i=a\Big)\leq \frac{1}{\sqrt{m(\pi)}}.
$$
Finally, we note that by Lemma~\ref{pi} we have many permutations with large $m(\pi)$, namely
$$
  \p_{\Pi_{2d}}\big(m(\pi)\leq k/50\big)\leq  \left( \frac{k}{1.1 d}\right)^{k/3}.
$$
Thus
$$
  \p\Big(\sum_{i=1}^d \xi_i=a\Big)\leq \sqrt\frac{50}{k} +
   \left( \frac{k}{1.1 d}\right)^{k/3}.
$$
This implies the desired result.
\qed

\medskip

\noindent{\bf Proof of Lemma~\ref{pi}. }
Without loss of generality we can assume that $x_i=1$ for $i\leq k$ and
$x_i=0$ for $i>k$. Let $\pi$ be a random permutation uniformly distributed
over $\Pi_{2d}$. The basic idea of the proof is to condition on a realization
of a set $\{i\leq d:\,x_{\pi(i)}=1\}$ and show that the conditional probability
of the event $|E|<k/50$ is small regardless of that realization.

Let $A=\{i\leq d:\,x_{\pi(i)}=1\}$ be a random subset of $[d]$.
Fix a subset $A_0\subset [d]$ with $|A_0|\leq k$
(then the event $A=A_0$ has a non-zero probability). Denote $m:=|A_0|$.
Further, define a random subset
$E_0=\{i\in A_0:\,x_{\pi(i+d)}=1\}$.
Clearly, we have $|E|=m-|E_0|+(k-m-|E_0|)=k-2|E_0|$ everywhere on the event $\{A=A_0\}$.
Let a parameter $\beta_1\in(0, 0.1)$ be chosen later.
Consider three cases.

\smallskip
\noindent
{\bf Case 1:  $m\leq(1-\beta_1)k/2$. \, }  Then clearly we have $|E|\geq k- 2m\geq \beta_1 k$
everywhere on the event $\{A=A_0\}$.

\smallskip
\noindent
{\bf Case 2: $m\geq(1+\beta_1)k/2$.\, } Since $|E_0|+m\leq k$ (deterministically),
we have $|E|\geq 2m-k\geq \beta_1 k$ everywhere on $\{A=A_0\}$.

\smallskip
\noindent
{\bf Case 3:   $(1-\beta_1)k/2\leq m\leq(1+\beta_1)k/2$.\, } Note that
the set $\{\pi(d+1),\pi(d+2),\dots,\pi(2d)\}$ contains $k-m$ ones and $d-k+m$ zeros. Thus, for every $\ell\leq k-m$ we have
$$
p_\ell:=\P(|E_0|=\ell\,|\,A=A_0)=\frac{1}{d!}{m\choose \ell}{d-m\choose k-m-\ell}(k-m)!(d-k+m)!,
$$
where the second factor is the number of choices of subsets $E_0$ of $A_0$ of cardinality $\ell$, the third factor is the number of
possible choices of the set $\{d+1\leq i\leq 2d:\,x_{\pi(i)}=1\mbox{ and }x_{\pi(i-d)}=0\}$ provided that $|E_0|=\ell$,
and the factors $(k-m)!$ and $(d-k+m)!$ are the numbers of all permutations of ones and zeros in the last $d$ positions.
Therefore,
$$
 p_\ell={d\choose m}^{-1}{k-m\choose \ell}{d-k+m\choose m-\ell}
$$
We choose $\beta >3/4$ from $(1-\beta)(1-\beta_1)/2=\beta _1$ and set $\beta_2=1-\beta$.
Using Chernoff bounds, we observe
\begin{align*}
\sum_{\ell\geq \beta m} {d-k+m\choose m-\ell} =
\sum_{\ell\leq (1-\beta) m} {d-k+m\choose \ell} \leq
\left(\frac{e(d-k+m)}{\beta_2 m}\right)^{\beta_2 m}
\leq \left(\frac{e d}{\beta_2 m}\right)^{\beta_2 m}.
\end{align*}
Since $k\leq 2m/(1-\beta_1)$ and $2/(1-\beta_1)-1-\beta=2\beta_2$,  then
for $\ell\geq \beta m$,
$$
 {k-m\choose \ell} = {k-m\choose k- m-\ell} \leq
 \left(\frac{e(k-m)}{k- m-\beta m}\right)^{k- m-\beta m}
 \leq \left(\frac{e(1+\beta_1)}{2(1-\beta_1)\beta_2}\right)^{2\beta_2 m}
 \leq \left(\frac{3}{2\beta_2}\right)^{2\beta_2 m}.
$$
Therefore we have
\begin{align*}
\sum_{\ell\geq \beta m}p_\ell\leq \left(\frac{m}{d}\right)^m \, \left(\frac{3}{2\beta_2}\right)^{2\beta_2 m} \, \left(\frac{e d}{\beta_2 m}\right)^{\beta_2 m} \leq \left(\frac{m}{d}\right)^{\beta m} \, \left(\frac{2}{\beta_2 }\right)^{3 \beta_2 m}\leq
\left(\frac{(1+\beta _1) k}{2d} \left(\frac{2}{\beta_2 }\right)^{3\beta_2/\beta}\right)^{\beta m}.
\end{align*}
Choosing $\beta _1=1/50$ we obtain
\begin{align*}
\sum_{\ell\geq \beta m}p_\ell\leq \left(\frac{k}{1.1 \, d}\right)^{k/3}.
\end{align*}
%
%
On the event $\{ A=A_0 \}$ we have $|E|\geq m-|E_0|$, hence
$$
 \P(|E|\leq (1-\beta)m\,|\,A=A_0)\leq\sum_{\ell\geq \beta m}p_\ell.
$$
Using that $m\geq(1-\beta_1)k/2$ and that $(1-\beta)(1-\beta_1)/2=\beta_1$ we obtain
$$
 \P(|E|\leq \beta_1 k\,\, |\, \,A=A_0)\leq
 \left(\frac{k}{1.1 \, d}\right)^{k/3},
$$
which completes the proof.
\qed

\subsection{Proof of the main theorem}
\label{main-result}

%
%
%

In this section, we complete the proof of the main result for adjacency matrices.
The general scheme is similar to the one in  \cite[Section~4]{Cook}.
The main idea of the proof of Theorem~A can be roughly described as follows:
after throwing away all small ``bad'' events (namely, the existence of almost
constant null vectors, big zero minors, and rows with largely overlapping supports)
we split the remaining singular matrices from $\Mc$ into two sets
$$
 E_1=\{M\in\Mc:\,\rk M= n-1\} \quad \mbox{ and } \quad E_2=\{M\in\Mc:\,\rk M\le n-2\}.
$$
Then, combining linear-algebraic arguments (Lemmas~\ref{l:n-2} and \ref{l:n-1}) with
the shuffling procedure (Lemma \ref{1-row}), we show that $E_1$ and $E_2$ have a small proportion inside the sets $\Mc$ and $\{M\in\Mc:\,\rk M\le n-1\}$, respectively. This implies that $E_1\cup E_2$ has small probability.

The argument is rather technical and involves various events and ``linear-algebraic''  objects. To make the reading more convenient, we first group the notation used in this section.

\subsubsection{Notation}

For every  $k\leq n$, let
$$
   \ek=\{M\in\Mc\,:\, \rk \M\leq k\}.
$$
Our purpose is to estimate the probability of the event $\en$.

Let $M$ be a matrix from $\Mc$ with rows $R_i$, $i\leq n$.
For every $i\in[n]$, we denote by $M^i$ the $(n-1)\times n$ minor of $M$
obtained by removing the row $R_{i}$.
Further, take a pair of distinct indices $(i, j)$, $i\ne j\le n$.
By $\MMij$ we denote the $(n-2)\times n$ minor of $M$ obtained by
removing the rows $R_i,R_j$. Additionally,  define
\begin{align*}
&\Vij=\Vij(M)=\spn\{R_k,\,k\neq i,j \}
\quad \mbox{ and } \quad
&\Fij=\Fij(M)=\spn\{\Vij,\,R_i+R_j \}.
\end{align*}
In what follows, we  write simply $\Vij$ and $\Fij$ instead of $\Vij(M)$ and $\Fij(M)$ as
the matrix $M$ will always be clear from the context.
Note  that  the random vector $R_i+R_j$ and the random subspaces $\Vij$ and $\Fij$
are fully determined by the $(n-2)\times n$ matrix $\MMij$.

As we  see later, to be able to successfully apply the aforementioned shuffling to a pair of rows $R_i,R_j$, we will need at our disposal a vector orthogonal to the subspace $\Fij$ such that its restriction to the union of the supports of $R_i$ and $R_j$ has many pairs of distinct coordinates. Of course, such a vector may not exist for some matrices $M\in\Mc$. We start by defining for every $q\in [n]$ and $i\ne j\le n$ a ``good'' subset of $\Mc$ as follows:
\begin{align}
\label{eijq}
\e^{i,j}(q)=\{M\in\Mc:\,\,&\exists v\perp \Fij
\quad\text{such that}
\\
&\forall\lam\in\R\,\,\, \, |\{k\in \supp (R_i+R_j):\;v_k\neq\lam\}|\geq q \}.\notag
\end{align}
For a matrix in this set
we fix one vector satisfying (\ref{eijq}), in fact we
define it as a function of the matrix. The crucial fact for our proof is
that since $\Fij$ and  $R_i+R_j$ are uniquely
determined by $\MMij$, we can fix such a vector for the
class of matrices ``sharing'' the same minor $\MMij$.
\begin{defi}\label{d:v}
 Given $M\in \e^{i,j}(q)$,   consider the equivalence class
$$
   {\cal{H}}_M^{i,j}(q) = \{\widetilde M\in \e^{i,j}(q) \, \, :\,  \, \widetilde M^{i,j} =M^{i,j} \}.
$$
For every equivalence class ${\cal{H}}_M^{i,j}(q)$  fix one vector
$v=v(M,q, i, j)$ satisfying
\begin{equation}
 v \perp F_{i,j} \;\; \mbox{and}
 \;\; \forall\lam\in\R\;\;|\{k\in \supp (R_i+R_j):\;v_k\neq\lam\}|\geq q.
\label{v}
\end{equation}
\end{defi}

Whenever $q$ and the indices $i,j$ are clear from context, we  write $v(M)$ instead of $v(M,q,i,j)$.
One of the key ideas of the proof of Theorem~A is to show
that for most matrices $M$ in ${\cal{H}}_M^{1,2}(q)$, the vector
$v(M)$ does not belong to the kernel of $M$. To this end we introduce a subset of $\ea$
\begin{equation*}
\K=\{M\in\ea:\, v(M)\in  \ker M\}.
\label{K}
\end{equation*}
In Lemma \ref{1-row} below we will show that the ratio $|\K|/|\Mc|$ goes to zero as  $d\to \infty$.

\medskip

As we mentioned above, in the proof we essentially restrict ourselves to the set of matrices,
which have no almost constant null-vectors,  no
big zero minors,  and no rows with largely overlapping supports. To define this set formally, let
$0< p\leq 1/3$, $2\leq q\leq d/2$, and
$\varepsilon\in (0, 1)$. Denote
\begin{align*}
\eg=\eg(p,q,\varepsilon):=\eac\cap{\Omega}^2_{\varepsilon}\setminus\eoo\big({p}/{2q}, {p}/{2}\big),
\label{eg}
\end{align*}
where  ${\Omega}^2_{\varepsilon}$,  $\eoo({p}/{2q}, {p}/{2})$, and $\eac$ were
introduced in  Proposition~\ref{lem-disjoint-rows},
\eqref{eo-def}, and \eqref{eAC}, respectively.
%
%
By Proposition~\ref{lem-disjoint-rows},  Theorem~\ref{t:AC} , and \eqref{remzero}
we have
\begin{equation}
     \p\big(\eg^c\big)\leq \frac{n^2}{2}\left(\frac{e d}{\eps n}\right)^{\eps d} +
      \left(\frac{C d}{n}\right)^{c d} + e^{-c n} \leq
      n^2\,  \left(\frac{e d}{\eps n}\right)^{\eps d}
      \leq  \left(\frac{e d}{\eps n}\right)^{\eps d/2}
       \label{nu}
\end{equation}
provided that $p\leq c_1/\ln d$, $q=c_2 p^2 d$, $c_3/\eps^2 \leq d\leq \eps n/6$ and $\eps$ is small enough.

\medskip

Further we will need two more auxiliary events dealing with the
$(n-2)\times n$ minors $\MMij$ of $M$. For $i\ne j,$ introduce
\begin{equation*}
\label{eijY}
 \e^{i,j}_{n-2}=\{M\in\Mc\,:\, \rk \MMij=n-2\quad\text{and}\quad R_i+R_j\notin V_{i,j}\},
\end{equation*}
and
\begin{equation*}
\label{rkij}
 \erkij=\{M\in\Mc\,:\, R_i, R_j\in V_{i,j}\}.
\end{equation*}
Note that  for every $M\in\erkij$ we have $\rk \M=\rk \MMij$. In the next section we prove several
statements involving events $\e^{i,j}_{n-2}$, $\erkij$, and $\K$.

\subsubsection{Proof of Theorem~A}

The next lemma encapsulates the shuffling procedure.
Recall that ${\Omega}^2_{\varepsilon}$ was defined in
Proposition~\ref{lem-disjoint-rows}.

\begin{lemma}\label{1-row}
Let $\varepsilon \in (0,1)$ and $2\varepsilon d< q\leq 2d/3$. Then
$$
  \p\left(\K\, \big| \,  \ea\cap{\Omega}^2_{\varepsilon}\right)\leq
  \frac{10}{\sqrt{(q-2\varepsilon d)}}.
$$
\end{lemma}

\smallskip

\begin{proof}
Note that
\begin{equation*}
 \K= \{M\in \ea:\, \langle v(M), R_1  \rangle =0 \}.
\end{equation*}
Let $M\in\ea\cap{\Omega}^2_{\varepsilon}$. Denote
\begin{align*}
  S_{1,2}=S_{1,2}(M)=\supp R_1\cup \supp R_2,\quad
  s_{1,2}=s_{1,2}(M)=\supp R_1\cap \supp R_2
\end{align*}
and set
$$
    S=S(M) = S_{1,2}\setminus s_{1,2},\quad   m_1= |S_{1,2}|, \quad m_2=|s_{1,2}|,
    \quad  \mbox{ and } \quad m = |S|.
$$
Note that $m_1=2d-m_2$ and $m=m_1-m_2 = 2(d-m_2)$.
By the definition of ${\Omega}^2_{\varepsilon}$, we have
$$
 m_1\geq 2(1-\varepsilon)d\quad \mbox{ and } \quad  m_2\leq 2\varepsilon d.
$$
By \eqref{v}, the vector $v:=v(M)$ satisfies $\, \forall \lam\in\R\quad |\{i\in S\,:\,x_i\ne\lam\}|\geq q-m_2$. Since $q-m_2\leq m/3$, by
Claim~\ref{c:J} (see Remark~\ref{aaaaaa})
there exists $J\subset S$ such that
\begin{equation}\label{vec-v}
  q-m_2\leq|J|\leq m-(q-m_2)  \;\text{ and }\; \forall i\in J\, \, \, \,
  \forall j\in S\setminus J \, \, \;\; v_i\neq v_j.
\end{equation}
 We compute the desired probability as follows.
For every (fixed) $M\in \ea\cap{\Omega}^2_{\varepsilon}$  consider the equivalence class
$$
   {\cal{F}}_M := {\cal{H}}_M^{1,2}(q) \cap{\Omega}^2_{\varepsilon}
= \left\{\widetilde M\in \ea\cap{\Omega}^2_{\varepsilon} \, \, :\,  \, \widetilde M^{1,2} =M^{1,2} \right\}.
$$
Note that by construction $S_{1,2}(\widetilde M)=S_{1,2}(M)=S_{1,2}$, $s_{1,2}(\widetilde M)=s_{1,2}(M)=s_{1,2}$ and $v(\widetilde M)=v(M)=v$
for every matrix $\widetilde M$ in ${\cal{F}}_M$. Therefore it is enough to show that the proportion of
matrices $\widetilde M \in {\cal{F}}_M$ satisfying
$\langle v, R_1(\widetilde M)\rangle =0$ is small.
Every matrix $\widetilde M \in {\cal{F}}_M$ is determined by its minor $M^{1,2}$ (which is
fixed on ${\cal{F}}_M$) and its first row $R_1(\widetilde M)$. Thus to determine a matrix in
${\cal{F}}_M$ it is enough to choose a support of the first row, which is a subset of $S_{1,2}$.
Since $m_2$ elements in $\supp R_1$ are fixed (as $s_{1,2}$ is fixed) then we have to calculate how many $(d-m_2)$-element subsets $B$ of $m$-element set $S$ exist so that
$$
    \langle v, R_1\rangle = \sum _{i\in B\cup s_{1,2}} v_i =0,
$$
that is
$$
  \sum _{i\in B} v_i =a := - \sum _{i\in s_{1,2}} v_i .
$$
For vectors $v=v(M)$ satisfying \eqref{vec-v} this was calculated in
Proposition~\ref{shary} (note that $m=2(d-m_2)$, $a$ is independent of
$B\subset S$ and apply the proposition with $q-m_2$ and $m/2$ instead of
$k$ and $d$). Applying this for all classes ${\cal{F}}_M$, we obtain the
desired bound.
\qed
\end{proof}

\smallskip

In what follows, we will show that, up to intersection with $\eg \cap \ea\cap \e^{1,2}_{n-2}$ (resp., $\eg \cap \ea\cap \erk$), the event
$E_1=\en\setminus \e_{n-2}$ (resp., $E_2= \e_{n-2}$) is a subset of $\K$, hence has a small probability. Our treatment of singular matrices $M$ with $\rk M=n-1$ and $\rk M\leq n-2$ is slightly  different, although the general idea is the same -- at the first step, given a singular matrix $M$, we fix a left null vector $y=y(M)$ and a right null-vector $x=x(M)$ and choose a row $R_i$ such that  $\rk M^i=\rk M$ and $x$ has many distinct coordinates on $R_i$. On the next step, we choose a second row $R_j$ so that the minor $M^{i, j}$ is of maximal rank, that is $\rk M^{i, j} = n-2$ in the case $\rk M=n-1$ and $\rk M^{i, j} = \rk M$ in the case $\rk M\leq n-2$. We also show that there are many choices for such $i$ and $j$. Finally, using the shuffling, we show, in a sense, that we can increase the rank of a matrix by ``playing" with rows $i$ and $j$ only, i.e. that the events $\e _{n-2}$ and $\en$ are small inside $\en$ and $\Mc$ respectively. The next lemma describes the set of ``good" $i$'s for the first step.


\medskip

\begin{lemma}
\label{l:1}
Let $0< p\leq 1/3$, and $ q\geq 2$ be such that $ pn/2q$ is an integer.
Further, let
$$
 M\in\en\cap\eac\setminus\eoo\big({p}/{2q}, {p}/{2}\big)
$$
and $x\in\ker M\setminus\{0\}$, $y\in\ker M^T\setminus\{0\}$.
Consider the set of indices
\begin{equation*}
 I_M(x,y)=\{i\in\supp y:\, \forall\lam\in\R\, \,
 |\{j\in\supp R_{i}\;:\;x_j\neq\lam\}|\geq q\}.
\end{equation*}
Then
\begin{equation*}
|I_M(x,y)|\geq pn/2.
 \label{|Ixy|}
  \end{equation*}
\end{lemma}

\medskip

 Note that for $y\in\ker M^T$ we have $\sum _i y_i R_i=0$ and
$I_M(x,y)\subset \supp y$. Therefore removing the row $R_i$,
$i\in I_M(x,y)$, we do not decrease the rank of $M$, that is
$\rk M^{i}=\rk M$.

\medskip

\noindent
{\bf Proof of Lemma \ref{l:1}. }
Since $x\notin AC(p)$ and $p\leq 1/3$, by Claim \ref{c:J} there
exists a subset $J_x\subset[n]$ such that
$$
   pn\leq |J_x|\leq (1-p)n \quad \text{and} \quad
   \forall i\in J_x \,\, \forall j\not\in J_x \quad x_i\neq x_j.
$$
Now we compute how many rows have more than $q$ ones in $J_x$ and
more than $q$ ones in $J_x^c$. Since $M\not\in \eoo({p}/{2q}, p/2)$
then applying Lemma~\ref{graph3} with $\alpha=p/(2q)$ and $\beta=p$,
we get
$$
  |\{i\;:\;|\supp R_i\cap J_x|\geq q\}|\geq(1- p/2q)n \quad \mbox{ and } \quad
  |\{i\;:\;|\supp R_i\cap J^c_x|\geq q\}|\geq(1- p/2q)n.
$$
Therefore
\begin{align*}
 |\{i\;:\;|\supp R_i\cap J_x|\geq q \quad \,\,\text{and}\quad \,\,|\supp R_i\cap
 J^c_x|\geq q\}|\geq(1- p/q)n.
\end{align*}
By the construction of the set $J_x$ this implies that the set
\begin{align*}
  I:=\{i\;:\;\forall\lam\in\R\quad|\{j\in \supp R_i\;:\;x_j\neq\lam\}|
   \geq q\}
\end{align*}
has cardinality at least $(1-p/q)n$.
Finally, since $y\notin AC(p)$, we have that $|\supp y|\geq pn$, which implies
$$
  |I_M(x,y)| = |I\cap\supp y| \geq pn -pn/q \geq pn/2,
$$
and completes the proof.
\qed

\bigskip

Now we consider the set of matrices $M\in \eg$ with $\rk M\leq n-2$.

\begin{lemma}\label{l:n-2}
Let $p,q$ satisfy the assumptions of Lemma~\ref{l:1}, $\eps\in (0, 1)$,
and let $\e=\enn\cap\eg$ with $\eg=\eg(p,q,\varepsilon)$. Then
\begin{equation*}
\p\big(\e\big)\leq2 p^{-2}\, \, \p\big(\erk\cap\e^{1,2}(q)\cap\e\big).
\end{equation*}
\end{lemma}

\smallskip

\noindent
{\bf Proof. }
Fix $M\in\e$. There exist $x\in\ker M\setminus\{0\}$ and
$y\in\ker M^T\setminus\{0\}$,  that is
$$
 \forall i\leq n \,\,\, x\perp R_i\quad \text{ and }\quad \sum_{i\in \supp y} y_i
 R_i=0.
$$
 Note that by the definition of $\eg$ we have $x,y\notin  AC(p)$. We compute
 how many ordered pairs $(i,j)$, $i\ne j$, satisfy
 \begin{equation*}
 R_i, R_j\in V_{i,j}\quad \text{ and }\quad \forall\lam\in\R\,\,|\{k\in \supp(R_i+R_j):\;x_k\neq\lam\}|\geq q.
   \end{equation*}
By Lemma~\ref{l:1}, the set $I_M(x,y)$ satisfies $|I_M(x,y)|\geq pn/2$,
and for every $i\in I_M(x,y)$ we have $\rk M^{i}=\rk M$.
Next, since $\rk M^{i}<n-1$, the set $\ker (M^{i})^T\setminus\{0\}$ is non-empty, i.e.\
$$
  \exists y^{(i)}\in\R^n\setminus\{0\}\, \, \mbox{ such that }\, \,  y^{(i)}_{i}=0,\quad \sum_{j=1}^n y^{(i)}_j R_j =0.
$$
Clearly $y^{(i)} \in\ker M^T\setminus\{0\}$, and, since $M\in \eac$,  $y^{(i)}$ has at least
 $pn$ non-zero coordinates; moreover,
$$
 \forall j\leq n \quad\text{ such that }\quad
  y^{(i)}_{j}\neq0\quad\text{ one has }\quad R_j\in V_{i,j}.
$$
This shows that for every $M\in \e$ there are at least $(pn)^2/4$ pairs
$(i, j)$ with $i<j$ satisfying $R_i, R_j\in V_{i,j}$. Obviously $x\perp\Fij$ for every $i,j\leq n$. Hence for each pair $(i, j)$ we have
\begin{equation*}
R_i, R_j\in V_{i,j}, \;\;  x\perp\Fij\,\, \, \mbox{ and } \, \, \, \forall\lam\in\R\,\,|\{k\in \supp(R_i+R_j):\;x_k\neq\lam\}|\geq q,
\end{equation*}
implying that $M$ belongs to at least $(pn)^2/4$ distinct subsets among $\{\erkij\cap\e^{i,j}(q)\cap\e\}_{i<j}$.

Thus, $\{\erkij\cap\e^{i,j}(q)\cap\e\}_{i<j}$
forms a ($(pn)^2/4$)-fold covering for $\e$. Since
$$
\p\left( \erkij \cap\e^{i,j}(q)\cap\e \right)
=
  \p\left( \erk \cap\e^{1,2}(q)\cap\e \right)
$$
 then applying Claim~\ref{k-fold}, we obtain
\begin{align*}
  \frac{(pn)^2}{4} \,  \p(\e)
&\leq
 \sum  _{i<j} \p\left( \erkij \cap\e^{i,j}(q)\cap\e \right)
=
  \frac{n(n-1)}{2} \, \p\left( \erk \cap\e^{1,2}(q)\cap\e \right).
\end{align*}
This  completes the proof of the lemma.

\qed

\medskip
Now we treat the case when a matrix has rank $n-1$.

\begin{lemma}\label{l:n-1}
Let $p,q$ satisfy the conditions of Lemma~\ref{l:1}, $\eps\in (0, 1)$,  and let
$\e=(\en\setminus\enn)\cap\eg$ with $\eg=\eg(p,q,\varepsilon)$. Then
\begin{equation*}
 \p\big(\e\big)\leq2 p^{-2}\, \, \p\big(\e^{1,2}_{n-2}\cap\e^{1,2}(q)\cap\e).
\end{equation*}
\end{lemma}

\noindent
{\bf Proof. }
Repeating the first step of the proof of Lemma \ref{l:n-2}, we fix $M\in\e$,
$x\in\ker M\setminus\{0\}$,  $ y\in\ker M^T\setminus\{0\}$. Then by
Lemma~\ref{l:1} the set of indices $I_M(x,y)$
has cardinality $|I_M(x,y)|\geq pn/2$ and for every $i\in I_M(x,y)$
the $(n-1)\times n$ minor $M^{i}$  satisfies
$\rk M^{i}=\rk M$.

Now, we calculate how many ordered pairs $(i,j)$, $i\ne j$,
exist such that
\begin{equation*}
  \rk M^{i,j}=n-2\quad \quad \mbox{ and } \quad \quad R_{i}+R_{j}\notin V_{i,j}.
\end{equation*}
Let  $i\in I_M(x,y)$. Since $y\notin AC(p)$, there are at least
$pn$ choices of $j$ such that $y_{j}\neq y_{i}$. Fix such a $j$. We claim that then
$R_{i}+R_{j}\notin V_{i,j}$.
Indeed, otherwise
$$
  R_{i}+R_{j}=\sum_{\ell\neq i,j}z_{\ell} R_{\ell}
$$
for some $z_\ell\in\R$,
hence there exists $w\in \ker M^T \setminus \{0\}$ such that $w_i=w_j$.
Since the dimension of $\ker M^T$ is one, we have $y=\lambda w$ for some $\lambda \in \R$,
which contradicts the condition $y_{i}\neq y_{j}$.
Therefore, there are at least $(pn)^2/2$ pairs $(i,j)$ with $i\ne j$ satisfying
$$
 \rk \MMij=n-2, \quad \quad R_{i}+R_{j}\notin V_{i,j},
$$
and
$$
  x\perp\Fij, \quad \quad \forall\lam\in\R\, \, \, \, |\{k\in \supp(R_{i}+R_{j})\;:\;x_k\neq\lam\}|\geq q.
$$
In other words, the matrix $M$ belongs to at least $(pn)^2/2$ events $\e^{i,j}_{n-2}\cap\e^{i,j}(q)\cap \e$.
Thus, we proved that
$\{\e^{i,j}_{n-2}\cap\e^{i,j}(q)\cap \e\}_{i<j}$
forms a $((pn)^2/4)$-fold covering of $\e$. Since for every $i<j$,
$$
  \p\left(\e^{i,j}_{n-2}\cap\e^{i,j}(q)\cap \e\right)=\p\left(\e^{1,2}_{n-2}\cap\e^{1,2}(q)\cap \e\right),
$$
then applying Claim~\ref{k-fold}, we obtain
\begin{align*}
\frac{(pn)^2}{4} \, \p\left( \e\right) \leq \sum_{i<j} \p\left(\e^{i,j}_{n-2}\cap\e^{i,j}(q)\cap \e\right)
= \frac{n(n-1)}{2} \, \p\left(\e^{1,2}_{n-2}\cap\e^{1,2}(q)\cap \e\right),
\end{align*}
and the proof is complete.
\qed

\medskip

We now  can finish Theorem~A.

\medskip

\noindent
{\bf Proof of Theorem~A. }
Let $p$, $q$, $\eps$ be chosen later to satisfy assumptions in the corresponding statements.
 By Lemmas~\ref{l:n-2}, \ref{l:n-1} and \eqref{nu} we obtain
\begin{align*}
  \p(\en)&\leq   \p\big(\enn\cap\eg\big)+ \p\big((\en\setminus \enn)\cap\eg\big)+ \left(\frac{e d}{\eps n}\right)^{\eps d/2}\\
  &\leq2p^{-2}\big(\p(A)+\p(B)\big)+ \left(\frac{e d}{\eps n}\right)^{\eps d/2},
 \end{align*}
 where
 \begin{align*}
  A=\erk\cap\e^{1,2}(q)\cap\enn\cap\eg \quad \mbox{ and } \quad
  B=\e^{1,2}_{n-2}\cap\e^{1,2}(q)\cap(\en\setminus \enn)\cap\eg.
 \end{align*}

We  show now  that
$$
A\cup B\subset \K\cap\e^{1,2}(q)\cap{\Omega}^2_{\varepsilon}.
$$
In other words, we verify that
for a matrix $M\in A\cup B$ the vector
$v(M)\in F^\perp_{1,2}$ (see Definition~\ref{d:v}) belongs to $\ker M$.

Indeed, if $M\in A$, then $R_1, R_2\in\V$.
Using the condition $v(M)\in F^\perp_{1,2}$ we immediately get $v(M)\in\ker M$.

 If $M\in B$ then  $\rk M=n-1$, $\rk M^{1,2}=n-2$, and $R_1+R_2\notin \V$.
 Therefore $\dim F_{1,2}=n-1$.  Since $\ker M\subseteq F^\perp_{1,2}$
and $\dim\ker M=\dim F^\perp_{1,2}=1$ we infer $\ker M= F^\perp_{1,2}$
and thus $v(M)\in\ker M$.

Finally note that $A$ and $B$ are disjoint, so $\p(A)+\p(B)=\p(A\cup B)$.
Applying Lemma \ref{1-row} we obtain
\begin{align*}
  \p(\en)\leq 2 p^{-2} \p\left(\K\cap\e^{1,2}(q)\cap{\Omega}^2_{\varepsilon}
  \right) + \left(\frac{e d}{\eps n}\right)^{\eps d/2}
  =\frac{20}{p^2\sqrt{(q-2\varepsilon d)}}+ \left(\frac{e d}{\eps n}\right)^{\eps d/2}.
\end{align*}
Finally we choose the parameters. Let $c_1$, $c_2$ be sufficiently small
positive absolute constants. Choose $p=c_1/\ln d$ and  $q$ to be the largest
integer not exceeding $c_2 p^2 d = c_1 c_2 d/\ln^2 d$ (we slightly adjust
$c_1, c_2$ so that $pn/2q$ is also an integer). Let $\eps = q/(4d)\approx 1/\ln ^2 d$
(note that then the condition $c_3/\eps^2 \leq d\leq \eps n/6$ needed in (\ref{nu}) is also satisfied).
Then we obtain the desired bound.
\qed

\begin{Remark}
We could choose $q=-c d p / \ln p\approx d /((\ln \ln d) \ln d)$, then
$\eps = q/(4d)\approx 1/((\ln \ln d) \ln d)$. This would lead to the restriction
$d\leq c n/((\ln \ln n) \ln n)$ instead of $d\leq c n/\ln^2 n$ in Theorem~A.
\end{Remark}

\address

\end{document}